\documentclass[PhD]{dukethesis}

\usepackage{graphicx}
\usepackage{amsmath}
\usepackage{amssymb}
\usepackage{amsthm}
\usepackage{amsfonts}

\newtheorem{prop}{Proposition}[section]
\newtheorem{cor}[prop]{Corollary}

\newtheorem{thm}[prop]{Theorem}
\newtheorem{lemma}[prop]{Lemma}
\newtheorem{conj}[prop]{Conjecture}

\theoremstyle{definition}
\newtheorem{defn}[prop]{Definition}

\newtheorem{ass}[prop]{Assumption}

\providecommand{\Label}[1]{\label{#1}}
\providecommand{\Ref}[1]{\ref{#1}}
\providecommand{\Cite}[1]{\cite{#1}}
\providecommand{\Eqref}[1]{\eqref{#1}}

\DeclareMathOperator{\Rc}{Rc}

\DeclareMathOperator{\Div}{div}

\DeclareMathOperator{\asinh}{asinh}
\DeclareMathOperator{\acos}{acos}

\newcommand{\fluff}[1]{\mbox{}}
\newcommand{\pp}[2]{\frac{\partial #1}{\partial #2}}
\newcommand{\dd}[2]{\frac{d#1}{d#2}}
\newcommand{\grad}{\nabla}
\newcommand{\p}{\partial}
\newcommand{\ga}{\alpha}
\newcommand{\gr}{\rho}
\newcommand{\gth}{\theta}
\newcommand{\gh}{\eta}
\newcommand{\gb}{\beta}
\newcommand{\gd}{\delta}
\newcommand{\gD}{\Delta}
\newcommand{\gf}{\varphi}
\newcommand{\gl}{\lambda}
\newcommand{\gs}{\sigma}
\newcommand{\gS}{\Sigma}

\newcommand{\gO}{\Omega}
\newcommand{\gv}{\nu}
\newcommand{\gt}{\tau}
\newcommand{\gP}{\Pi}
\newcommand{\gp}{\pi}
\newcommand{\gPs}{\Psi}
\newcommand{\gps}{\psi}
\newcommand{\gc}{\chi}
\newcommand{\gk}{\kappa}
\newcommand{\gm}{\mu}

\newcommand{\gz}{\zeta}

\newcommand{\N}{\nabla}
\renewcommand{\ge}{\epsilon}
\renewcommand{\gg}{\gamma}
\renewcommand{\Re}{\ensuremath{\mathbb{R}}}
\newcommand{\Ree}{\mathbb{R}\text{e}}
\newcommand{\Iem}{\mathbb{I}\text{m}}
\newcommand{\Ce}{\ensuremath{\mathbb{C}}}
\newcommand{\abs}[1]{\ensuremath{\left|{#1}\right|}}
\newcommand{\norm}[1]{\ensuremath{\left\|{#1}\right\|}}

\newcommand{\ADM}{\text{ADM}}
\newcommand{\Hk}{\text{H}}
\newcommand{\B}{\text{B}}
\newcommand{\NPMS}{\text{R}}

\newcommand{\E}{\mathcal{E}}
\newcommand{\de}[1]{\,d{#1}}
\renewcommand{\bar}{\overline}
\renewcommand{\tilde}{\widetilde}
\newcommand{\bigo}{\mathcal{O}}
\newcommand{\ip}[2]{\left<{#1},{#2}\right>}

\author{Nicholas P. Robbins}

\title{Negative Point Mass Singularities in General Relativity}

\advisor{Hubert L. Bray}

\department{Mathematics}

\subject{Mathematics}

\date{2007}

\member{Mark A. Stern}
\member{Paul A. Aspinwall}
\member{N. Ronen Plesser}

\begin{document}


\maketitle

\makeabstract
\Copyright

\abstract

First we review the definition of a negative point mass
singularity. Then we examine the gravitational lensing effects of
these singularities in isolation and with shear and convergence from
continuous matter.  We review the Inverse Mean Curvature Flow and use
this flow to prove some new results about the mass of a singularity,
the ADM mass of the manifold, and the capacity of the singularity. We
describe some particular examples of these singularities that exhibit
additional symmetries.

\acknowledgements

I want to thank my thesis advisor, Professor Hubert Bray.  Without
your copious support, advice, and encouragement throughout my thesis
process this would not have been possible.

Thanks also to the many faculty members who have helped me throughout
my career at Duke.  Thanks to Professors Mark Stern, Paul Aspinwall
and Ronen Plesser for your extensive feedback during the thesis
writing process. I want to thank Professor Stern in particular for the
courses which formed the backbone of my coursework at Duke.  Thanks to
Professor Tom Witelski for your helpful advice and support at a
difficult time in my graduate career.

I have had wonderful support for developing my teaching skills while
at Duke. Thanks to James Tomberg and Professor Jack Bookman in
particular for your extensive support for my teaching and career, as
well as Professor Lewis Blake for your work supporting undergraduate
instruction.

I also want to thank Professor Jennifer Hontz at Meredith College for
allowing me to see academic life outside of Duke University.

I want to thank all the departmental staff. In particular, thanks to
Georgia Barns and Shannon Holder for helping me over innumerable
institutional hurdles.

I have had the pleasure of many meaningful friendships with the other
graduate students in the program. In particular I wish to thank my
officemates Melanie Bain, Thomas Laurent, Michael Nicholas, Michael
Gratton, and Ryan Haskett, as well as Jeffrey Streets, Abraham Smith,
Paul Bendich, and Greg Firestone for many insightful conversations
academic and otherwise. Without your friendships the road these six
years would have felt far rougher.

I have the good fortune of coming from a supportive and loving family.
Without the support of my mother and sister, I could not have even
gotten to graduate school, much less finished.

Finally, I wish to thank my partner, Claire. I can't imagine how
difficult these years would have been without your ever-present faith,
encouragement, support and humor.

\tableofcontents

\listoftables	

\listoffigures	

%
%
%
\chapters

\chapter{Introduction}

General Relativity is currently the accepted  model of the physics of the
universe on large scales. The theory principally consists of three parts. 
\begin{itemize}
\item The universe is modeled by a four dimensional Lorentzian
  manifold. 
\item The geometry of this manifold is given by the Einstein equation:
  \begin{gather}
    \Rc_{\gm\gv}-\frac12Rg_{\gm\gv}=8\pi T_{\gm\gv}.
  \end{gather}
  The left hand side is called the Einstein Tensor and the right hand
  side is the stress energy tensor for the system. The stress energy
  tensor is a property of the matter in the system.
\item Particles, in the absence of other forces, move on timelike
  geodesics. 
\end{itemize}
The theory has had remarkable success. Its major accomplishments
include explaining the orbital perihelion precession of Mercury,
differing values for the bending of light near a massive body, and a
working cosmological model running up to moments after the Big Bang.

In particular, Gravitational Lensing has provided the most useful
method for observing and studying that portion of the matter in the
universe that does not emit or absorb light. This is currently
believed to comprise approximately 85\% of the total matter in the
universe.

Despite its many successes, there are a number of difficulties that
arise when working in the theory. The Einstein Equation, even in a
vacuum, is a nonlinear second order hyperbolic differential equation
in the metric. This makes exact solutions difficult to find. While the
behavior of small test particles is easily given by the geodesic
condition, the behavior of continuous masses, with features like
tension, pressure, et cetera, are not given by the theory, but require
an external derivation. Furthermore the relevant manifold is Lorentzian,
not Riemannian, removing many powerful tools.

To avoid many of these difficulties, one may study Riemannian General
Relativity. This is the study of Riemannian 3-manifolds that could
arise as spacelike hypersurfaces in a spacetime in General Relativity.
The Einstein Equation is translated into equations about the metric on
this spacelike hypersurface and its second fundamental form. The
properties of the matter in the theory are replaced by conditions such
as the Dominant and Weak Energy conditions, which are also conditions
on the metric of the hypersurface and its second fundamental form.
Many of the important questions in General Relativity have analogues
on Riemannian General Relativity. For example, the Penrose Conjecture
can be restricted to the Riemannian Penrose Inequality.

This thesis consists of a study of the properties of Negative Point
Mass Singularities. The motivating example of which is the spacial
Schwarzschild metric with a negative mass parameter: 
\begin{gather}
 g_{ij} = \left(1+\frac{m}{2r}\right)^4\gd_{ij}\qquad\qquad m<0.
\end{gather}
In addition to being historically and physically important, the
Schwarzschild solution is of particular mathematical interest since it
is the case of equality of the Riemannian Penrose conjecture, and, in
the case when $m=0$, it is the case of equality of the Riemannian
Positive Mass Theorem. Thus this metric, and its generalizations,
show promise as objects of study.

This thesis consists of three main topics. After laying out the
necessary definitions in Chapter~\Ref{Defns}, we examine the
gravitational lensing effects of these singularities in
Chapter~\Ref{Lens}. Next we summarize the work of Huisken and Ilmanen
on Inverse Mean Curvature Flow in Chapter~\Ref{IMCF}. We then use this
to prove a number of results in Chapter~\Ref{IMCF-NPMS}. In
Chapter~\Ref{Axi} we show what additional information can be gained if
the singularities possess additional symmetries.

\chapter{Definitions} 
\Label{Defns}
\section{Asymptotically Flat Manifolds}

Physically, we want to make sure our manifolds represent isolated
systems. A precise formulation of this is asymptotic flatness.
\begin{defn}(\Cite{IMCF})
  A Riemannian 3-manifold $(M,g)$ is called \emph{asymptotically flat}
  if it is the union of a compact set $K$, and sets $E_i$
  diffeomorphic to the complement of a compact set $K_i$ in $\Re^3$,
  where the metric of each $E_i$ satisfies
\begin{gather}
  \abs{g_{ij}-\gd_{ij}}\leq \frac{C}{\abs{x}},\qquad
  \abs{g_{ij,k}}\leq \frac{C}{\abs{x}^2}
\end{gather}
as $\abs{x}\to\infty$. Derivatives are taken in the flat metric
$\gd_{ij}$ on $x\in\Re^3$. Furthermore the Ricci curvature must
satisfy 
\begin{gather}
  \Rc \geq -\frac{Cg}{\abs{x}^2}.
\end{gather}
The set $E_i$ is called an \emph{end} of $M$.
\end{defn}
A manifold may have several ends, but most of our results will be
relative to a single end.

In \cite{adm}, Arnowitt, Deser and Misner define a geometric invariant
that is now called the ADM mass.
\begin{defn}
The \emph{ADM mass} of an end of an asymptotically flat manifold is 
\begin{gather}
  m_{\ADM}=\lim_{r\to\infty}\frac1{16\pi}\int_{S_\gd^r}\left(g_{ij,i}
    -g_{ii,j}\right)n^j d\gm.
\end{gather}
\end{defn}
This quantity is finite exactly when the total scalar curvature of the
chosen end is finite. This definition appears to be coordinate
dependent, however in \Cite{adm} the authors show that it is actually
an invariant when
\begin{gather}
  \int_{M\setminus K}\abs{R}<\infty.
\end{gather}
\section{Quasilocal Mass Functionals}
While the ADM mass provides a definition for the total mass of a
manifold, or the mass seen at infinity, there is no computable
definition for the mass of a region. The two of most relevance are the
(Riemannian) Hawking mass and the Bartnik mass. 
\begin{defn}
  The \emph{Hawking mass} of a surface $\gS$ is given by
  \begin{gather}
    m_{\Hk}=\sqrt{\frac{\abs{\gS}}{16\pi}}\left(1-\frac1{16\pi}\int_\gS
      H^2\right).
  \end{gather}
\end{defn}
Consider the Hawking mass of a surface, $\gS$, in $\Re^3$, and note
that $H=\gk_1+\gk_2$, where $\gk_i$ are the principal curvatures of
$\gS$. Thus $H^2 = \gk_1^2+\gk^2_2 + 2\gk_1\gk_2$. Hence, for a sphere in $\Re^3$,
\[
\int_\gS
H^2=\int_\gS 2K +\gk_1^2+\gk^2_2\geq\int_\gS 2K + 2K =
8\pi\gc(\gS)= 16\pi.
\]
Thus the Hawking mass of a sphere is always nonpositive in $\Re^{3}$.
Furthermore the Hawking mass can be decreased by making the surface
$\gS$ have high frequency oscillations. These two observations lead to
the conclusion that the Hawking mass tends to underestimate the mass
in a region.

The other quasilocal mass functional of interest is the
Bartnik mass defined in \Cite{bartnik-mass}.
\begin{defn}
  Let the asymptotically flat manifold $(M,g)$ have nonnegative scalar
  curvature.  Let $\gO$ be a domain in $M$ with connected boundary.
  Assume $M$ has no horizons (minimal spheres) outside of $\gO$. Call
  any asymptotically flat manifold $(\tilde M, \tilde g)$ acceptable
  if it has nonnegative scalar curvature, contains an isometric copy
  of $\gO$, and has no horizons outside of $\gO$.  Then the
  \emph{Bartnik mass}, $m_{\B}(\gO)$ of $\gO$ is defined to be the
  infimum of the ADM masses of these acceptable manifolds.
\end{defn}
The positive mass theorem guarantees that this mass will be positive
if the interior of $\gS$ fulfills the hypotheses of the theorem.  The
Bartnik mass is very difficult to compute. The only cases where it is
known are when the surface can be embedded in the exterior region of
the Schwarzschild spacial metric or in $\Re^3$.  The Schwarzschild
metric is the case of equality of the Riemannian Penrose inequality
and $\Re^3$ is the case of equality for the positive mass theorem.


\section{Definition and Mass of Negative Point Mass Singularities}

The basic example of a negative point mass singularity is the negative
Schwarzschild solution. This is the manifold $\Re^3\setminus B_{-m/2}$
with the metric
\begin{gather}
  g_{ij}  = \left(1+\frac{m}{2r}\right)^{4}\gd_{ij}
\end{gather}
where $m<0$. This manifold fails the requirements of the positive
mass theorem since it is not complete: geodesics reach the sphere at
$r=-m/2$ in finite distance. A straightforward calculation shows that
the ADM mass of this manifold is given by $m$.  Furthermore the far
field deflection of geodesics is the same as for a Newtonian mass of
$m$.  These results are identical to the same results for a positive
mass Schwarzschild solution. 

Two important aspects of this example will be incorporated into the
definition of a negative point mass singularity. One is that the point
itself is not included. To justify the use of the word ``point'' we
have to describe the behavior of surfaces near the singularity. The
manifold in that region should have surfaces whose areas converge to
zero. In addition the capacity of these surfaces should go to zero.
The second aspect is the presence of a background metric, in this case
the flat metric.  This background metric will provide a location where
we can compute information about the singularity.

This example motivates the following definition.
\begin{defn}
  Let $M^3$ be a smooth manifold with boundary, where the boundary is
  compact.  Let $\gP$ be a compact connected component of the boundary
  of $M$. Let the interior of $M$ be a Riemannian manifold with smooth
  metric $g$. Suppose that, for any smooth family of surfaces which
  locally foliate a neighborhood of $\gP$, the areas with respect to
  $g$ go to zero as the surfaces converge to $\gP$.  Then $\gP$ is a
  \emph{Negative Point Mass Singularity.}
\end{defn}
A particularly useful class of these singularities are Regular
Negative Point Mass Singularities.
\begin{defn}
  Let $M^3$ be a smooth manifold with boundary. Let the boundary of
  $M$ consist of one compact component, $\gP$.  Let $(M^3\setminus
  \gP,g)$ be a smooth Riemannian manifold.  Suppose that, for any
  smooth family of surfaces which locally foliate a neighborhood of
  $\gP$, the areas with respect to $g$ go to zero as the surfaces
  converge to $\gP$.  If there is a smooth metric $\bar g$ on $M^3$
  and a smooth function $\bar \gf$ on $M$ with nonzero differential on
  $\gP$ so that $g=\bar\gf^4\bar g$, then we call $\gP$ a
  \emph{Regular Negative Point Mass Singularity.}  We call the data $(
  M^3,\bar g, \bar \gf)$ a \emph{resolution} of $\gP$.
\end{defn}
Notice that while $\gP$ is topologically a surface, and it is a
surface in the Riemannian manifold $(M^3,\bar g)$, the areas of
surfaces near it in $(M^3\setminus\gP,g)$ approach zero, so we will
sometimes speak of $\gP$ as being a point $p$, when we are thinking in
terms of the metric $g$. Furthermore, notice that the requirement that
areas near $\gP$ go to zero under $g$ tells us that $\bar\gf=0$ on
$\gP$.

We can define the mass of a Regular Negative Point Mass Singularity as follows:
\begin{defn}
  Let $(M^3,\bar g,\bar \gf)$ be a resolution of a regular negative
  point mass singularity $p=\gP$. Let $\bar\gv$ be the unit normal to
  $\gP$ in $\bar g$.  If the capacity of $p$ is zero, then the
  \emph{regular mass} of $p$ is defined to be
  \begin{gather}
    m_{\NPMS}(p)=-\frac14\left(\frac1\pi\int_{\gP}\bar\gv(\bar
      \gf)^{4/3}\, \bar{dA}\right)^{3/2}.
  \end{gather}
  If the capacity of $p$ is nonzero, then the mass of $p$ is defined to be
  $-\infty$.
\end{defn}
See Chapter \ref{IMCF-NPMS} for a discussion of the capacity of
points like $p$.  We can also define the mass of a negative point mass
singularity that may not be regular.
\begin{defn}\Label{NPMS-mass-surf}
  Let $(M^3,g)$ be an asymptotically flat manifold, with a negative
  point mass singularity $p$.  Let $\gS_i$ be a smooth family of
  surfaces converging to $p$. Define
  $h_i$ by
  \begin{align}
    \gD h_i &= 0 \\
    \lim_{x\to\infty}h_i &=1\\
    h_i &= 0 \text{ on } \gS_i .
  \end{align}
  Then the manifold $(M,h^4_ig)$ has a negative point mass singularity
  at $\gS_i=p_i$ which is resolved by $(M,g,h_i)$.  Define the mass of
  $p$ to be
  \begin{gather}
    \sup_{\{\gS_i\}}\varlimsup_{i\to \infty}
    -\frac14\left(\frac1\pi\int_{\gS_i} \gv(h_i)^{4/3}\, dA\right)^{3/2}
    =\sup_{\{\gS_i\}}\varlimsup_{i\to \infty} m_{\NPMS}(p_i).
  \end{gather}
  Here the outer sup is over all possible smooth families of surfaces
  $\{\gS_i\}$ which converge to $p$. 
\end{defn}
A straightforward calculation shows that if the capacity of $p$ is
non-zero, then the mass of $p$ is $-\infty$.  This is the definition
we will be working with. However, an alternative definition for the
mass of a negative point mass singularity follows.

\begin{defn}
  Let $(M^3,g)$ be an asymptotically flat manifold, with a negative
  point mass singularity $p$. Choose a function $h$ that satisfies the
  equations
  \begin{align}
    \gD h &= 0 \\
    h &=\frac{1}{r}+\mathcal O\left(\frac{1}{r^2}\right)\\
    \lim_{x\to p}h &= \infty.
  \end{align}
  Now define surfaces $\gS_t=\left\{x|h(x)=t\right\}$ and functions
  $\gf_t(x)=1-h/t$. 
  Then the manifold $(M,\gf^4_tg)$ has a negative point mass singularity
  at $\gS_t=p_t$ which is resolved by $(M,g,\gf_t)$.  Define the mass of
  $p$ to be
  \begin{gather}
    \sup_{h}\varlimsup_{t\to \infty}
    -\frac14\left(\frac1\pi\int_{\gS_t} \gv(\gf_t)^{4/3}\, dA\right)^{3/2}
    =\sup_{h}\varlimsup_{t\to \infty} m_{\NPMS}(p_t).
  \end{gather}
  Here the outer sup is over all possible $h$'s.
\end{defn}

\section{Fundamental Results}

Before we continue we must verify that these definitions are
consistent.  First it must be verified that the regular mass of a
regular singularity is indeed intrinsic to the singularity, as shown
in \Cite{Hugh-INI}.
\begin{lemma}
  The regular mass of a negative point mass singularity is
  independent of the resolution. 
  \begin{proof}
    Let $(M^3,\bar g, \bar \gf)$ and $(M^3,\tilde g, \tilde \gf)$ be
    two resolutions of the same negative point mass singularity, $p$. 
    Then define $\gl$ by
    $\bar \gf = \gl \tilde \gf$. Thus we note the following scalings:
    \begin{gather}
      \tilde g = \gl^4 \bar g\\
      \tilde{dA}= \gl^4\bar{dA}\\
      \tilde \gf = \gl^{-1}\bar\gf\\
      \tilde \gv = \gl^{-2}\bar\gv
    \end{gather}
    Now note that since $\tilde\gf,\bar\gf=0$ on $\tilde\gP,\bar\gP$, 
    \begin{gather}
      \tilde{\gv}(\tilde\gf)= \gl^{-2}\bar\gv\left(\gl^{-1}\bar\gf\right)=
      \gl^{-3}\bar\gv\left(\bar\gf\right)+\gl^{-4}\bar\gv(\gl)\bar\gf. 
    \end{gather}
    The last term, $\gl^{-4}\bar\gv(\gl)\bar\gf$, needs discussion.
    Both $\bar\gf$ and $\tilde\gf$ are smooth functions with zero set
    $\gP$ and they both have nonzero differential on $\gP$. Thus $\gl$
    is smooth. To see this choose a coordinate patch on the boundary
    where $\gP$ is given by $x=0$. Then Taylor's formula tells us that
    \begin{gather}
      \gl = \frac{\int_0^1 \pp{\bar\gf}{x}(xs,y,z)\,ds}{\int_0^1
        \pp{\tilde\gf}{x}(xs,y,z)\,ds},
    \end{gather}
    which is a nonzero smooth function. Thus since $\bar\gf$ goes to
    zero on $\gP$, this last term is zero on $\gP$.  Thus the mass of
    $p$ using the $(M^3,\tilde g,\tilde\gf)$ resolution is given by
    \begin{align}
      m_{\NPMS}(p)&=-\frac14\left(\frac1\pi\int_{\tilde\gP}\tilde\gv(\tilde
        \gf)^{4/3}\, \tilde{dA}\right)^{3/2}\\
      &=-\frac14\left(\frac1\pi\int_{\bar\gP}\left[\gl^{-2}
        \bar\gv(\gl^{-1}\bar\gf)\right]^{4/3}\, \gl^4\bar{dA}\right)^{3/2}\\
      &=-\frac14\left(\frac1\pi\int_{\bar\gP}\left[\gl^{-3}
        \bar\gv(\bar\gf)\right]^{4/3}\, \gl^4\bar{dA}\right)^{3/2}\\
      &=-\frac14\left(\frac1\pi\int_{\bar\gP}
        \bar\gv(\bar\gf)^{4/3}\,\bar{dA}\right)^{3/2}.
  \end{align}
  \end{proof}
\end{lemma}

Definition~\Ref{NPMS-mass-surf} seems to involve the entire manifold,
as the definition of $h_i$ takes place on the entire manifold. However
that isn't the case. The mass is actually local to the point $p$.
\begin{lemma}
  Let $(M^3,g)$ be a manifold with a negative point mass singularity
  $p$. Let $\tilde g$ be a second metric on $M$ that agrees with $g$
  in a neighborhood of $p$. Then the mass of $p$ in $(M^3,g)$ and
  $(M^3,\tilde g)$ are equal.
  \begin{proof}
    The goal is to show that for any selection of
    $\left\{\gS_i\right\}$, the series $ m_{\NPMS}(p_i)$ and $\tilde
    m_{\NPMS}(p_i)$ obtained in the calculation of the mass of $p$,
    with respect to $(M,g)$ and $(M,\tilde g)$ converge.  Let $S$ be a
    smooth, compact, connected surface separating $p$ from infinity
    and contained in the region where $g$ and $\tilde g$ agree.  Fix
    $i$ large enough so that $\gS_i$ is inside of $S$, and suppress
    the index $i$ on all our functions.  Then define the functions
    $h,\tilde h$ by
    \begin{align*}
      h=\tilde h &= 0\text{ on } \gS_i\\
      \lim_{x\to\infty}h=\lim_{x\to\infty}\tilde h & = 1\\
      \gD h = \tilde \gD \tilde h & = 0.
    \end{align*}
    Here $\gD$ and $\tilde \gD$ denote the Laplacian with respect to
    $g$ and $\tilde g$ respectively.

    Now inside $S$, $\gD=\tilde \gD$ since $g=\tilde g$. Thus there is
    only one notion of harmonic, and $h$ and $\tilde h$ differ only by
    their boundary values on $S$.  Let $\ge=1-\min_{S}\{h,\tilde h\}$.
    Consider the following two functions $f^-$ and $f^+$ defined between
    $S$ and $\gS_i$:
    \begin{align*}
      f^-=f^+&=0\text{ on }\gS_i\\
      \gD f^- = \gD f^+ &=0\\
      f^-&=1-\ge\text{ on } S\\
      f^+&=1\text{ on } S.
    \end{align*}
    Thus by the maximum principle, we have the following inside $S$
    \begin{gather}
      f^+\geq h,\tilde h\geq f^-.
    \end{gather}
    Furthermore, since all four functions are zero on $\gS_i$,
    \begin{gather}
      \gv(f^+)\geq \gv(h),\gv(\tilde h)\geq \gv(f^-).
    \end{gather}
    Here $\gv$ is the normal derivative on $\gS_i$. Now define
    $\E(\gf)$ by the formula
    \begin{gather}
      \E(\gf) = \int_{\gS_i} \gv(\gf)^{4/3}\, dA.
    \end{gather}
    Then the ordering of the derivatives gives the ordering
     \begin{gather}
      \E(f^+)\geq \E(h),\E(\tilde h)\geq \E(f^-).\Label{nrg-sand}
    \end{gather}
    However, since $f^-= (1-\ge)f^+$, 
    \begin{gather}
      \gv(f^-) = (1-\ge)\gv(f^+),
    \end{gather}
    hence,
    \begin{gather}
      \E(f^-) = (1-\ge)^{4/3} \E(f^+).
    \end{gather}
    Now, without loss of generality assume that the limit of the
    capacities of $\left\{\gS_i\right\}$ is zero, as the mass would be
    $-\infty$ otherwise.  Thus as $i\to\infty$, $\gS_i$ has
    capacity going to zero. Hence $\ge_i$ goes to zero, and so $
    \E(f^-_{i})/\E(f^+_i)$ goes to 1. Thus equation~\Eqref{nrg-sand} forces
    $\E(h_i)$ and $\E(\tilde h_i)$ to equality.  This forces the
    masses of $p_i$ in the two metrics to equality as well.
  \end{proof}
\end{lemma}
\begin{cor}
  In Definition~\Ref{NPMS-mass-surf} we may replace the condition that
  $\gf_i$ be one at infinity with the condition that $\gf_i$ be one
  on a fixed surface outside $\gS_i$ for $i$ sufficiently large.
\end{cor}

This mass also agrees with the regular mass when the singularity is
regular.

\begin{lemma}[\Cite{Bray-unpub}]
  Let $(M^3,g)$ be an asymptotically flat manifold with negative point
  mass singularity $p$. Let $p$ have a resolution $(M^3,\tilde g,
  \tilde\gf)$. Then the regular mass of $p$ equals the general mass of
  $p$.
\fluff{
  \begin{proof}
    We can assume that the capacities of the $\gS_i$ go to zero, since
    otherwise both masses are $-\infty$. Assume $\gf$ goes to 1 at
    infinity, and is bounded with bound $M$. Let $\{\gS_i\}$ be a
    sequence of surfaces converging to $p$. Let $h_i$ be the harmonic
    function going to 1 at infinity and 0 on $\gS_i$.  Extend $h_i$ to
    $p_i$ by zero.  Since $\gS_i$ is outside $\gS_{i+1}$, the maximum
    principle tells us that $h_i\leq h_{i+1}$. So the $h_i$ are
    nondecreasing in $i$.  Since the capacity of $p$ is zero, the
    $h_i$ are converging pointwise to $1$.  Now we define the
    functions $\tilde\gf_i=h_i\tilde\gf$.  These functions are smooth
    off of $\gS_i$ since $\tilde\gf$ and $h_i$ are. They converge
    pointwise to $\tilde\gf$.  To see that they converge uniformly to
    $\tilde\gf$, let $\ge>0$ be chosen. Choose $I$ so that
    $\tilde\gf<\ge$ inside $\gS_i$ for all $i>I$. Choose $J$ so that
    $1-h_j<\ge/M$ outside of $\gS_I$ for all $j>J$.  Then, for all
    $k>J,I$, outside of $\gS_I$, we have
    \begin{gather}
      \abs{\tilde\gf-\tilde\gf_k}\leq \frac{\ge}{M}\abs{\tilde\gf}\leq \ge,
    \end{gather}
    while inside $\gS_I$ we have
    \begin{gather}
      \abs{\tilde\gf-\tilde\gf_k}= (1-h_k)\tilde\gf <\ge.
    \end{gather}
    Since $\tilde\gf_i\to\tilde\gf$ uniformly
    \begin{gather}
      \lim_{j\to\infty} \int_{\gS_i} \tilde\gv(\tilde\gf_j)^{4/3}\, dA
      = \int_{\gS_i} \tilde\gv(\tilde\gf)^{4/3}\, dA
    \end{gather}
    uniformly in $i$. 
    Since $\tilde\gf$ is smooth
    \begin{gather}
      \lim_{i\to\infty}\int_{\gS_i} \tilde\gv(\tilde\gf)^{4/3}\, dA
      = \int_{\gS} \tilde\gv(\tilde\gf)^{4/3}\, dA.
    \end{gather}
    Thus the integrals
    \begin{gather}
      \int_{\gS_i}\tilde\gv(\tilde\gf_i)^{4/3}\,dA
    \end{gather}
    converge to the integral
    \begin{gather}
      \int_{\gS} \tilde\gv(\tilde\gf)^{4/3}\, dA
    \end{gather}
    as desired.

   \fluff{ also satisfy the equations
    \begin{gather}
      \begin{aligned}
        \tilde\gD\tilde\gf_i &= h_i\tilde\gD\tilde\gf\\
        \lim_{x\to\infty}\tilde\gf_i &= 1\\
        \tilde\gf_i(\gS_i)&=0.
      \end{aligned}
      \Label{pde-loc}
    \end{gather}
So the $\tilde\gf_i$ are bounded and, outside of $\gS_i$, By the
    maximum principle they are nondecreasing.  Thus the $\tilde\gf_i$
    converge uniformly.  Since as $i$ goes to infinity, $\gS_i\to\gS$,
    Equation~\Eqref{pde-loc} becomes the same as the defining equation
    for $\tilde\gf$. Thus $\lim_{i\to\infty}\tilde\gf_i=\tilde\gf$.
    Thus
    \begin{gather}
      \lim_{j\to\infty} \int_{\gS_i} \tilde\gv(\tilde\gf_j)^{4/3}\, dA
      = \int_{\gS_i} \tilde\gv(\tilde\gf)^{4/3}\, dA.
    \end{gather}
    Since $\tilde\gf$ is the solution to an elliptic equation, its
    gradient is bounded. Thus since as $i\to\infty$ $\gS_i$ goes to $\gS$,
    \begin{gather}
      \lim_{i\to\infty}\int_{\gS_i} \tilde\gv(\tilde\gf)^{4/3}\, dA
      = \int_{\gS} \tilde\gv(\tilde\gf)^{4/3}\, dA,
    \end{gather}
    and our two masses are equal. 
}
\end{proof}
  }
\end{lemma}


\chapter{Gravitational Lensing by Negative Point Mass Singularities}
\Label{Lens}

\section{Gravitational Lensing Background}
One of the first testable predictions of general relativity was the
difference in the deflection of light by gravity. This effect was
first confirmed during the 1919 solar eclipse. Since then
gravitational lensing has become an powerful tool for astronomy in
general and cosmology in particular. Gravitational lensing has made it
possible to detect the presence of dark matter by observing its
effects on background images.

In this chapter we will develop the properties of lensing by negative
point mass singularities in the setting of accepted cosmology.  We
will make a number of assumptions based on that cosmology that will
allow us to obtain simple formulas for gravitational lensing by
negative point mass singularities. Then we will characterize their
lensing effects and compare them to positive mass point sources.  We
will not cover the entire field of gravitational lensing, but only
develop enough for our purposes. 

We will restrict ourselves to negative point mass singularities that
agree with a negative Schwarzschild solution to first order. We will
find that the lensing effects of these singularities in the presence
of continuous matter and shear can be duplicated by configurations
with positive mass lenses.

We will follow the presentation given in \Cite{plw}. We will differ
from this presentation by using geometrized units where $c=G=1$. We
will also consider lens potentials outside of the scope of \Cite{plw}.

\section{Cosmology for Gravitational Lensing} 
To simplify the calculations involved in lensing, we will make use of
a number of assumptions about the configuration of our system and its
behavior. These assumptions are based on the scales and phenomenology
of astronomy. Our first assumption is one of cosmology.
\begin{ass}
  The universe is described by a Friedmann--Robertson--Walker cosmology.
\end{ass} 
The Friedmann--Robertson--Walker model is an isotropic homogeneous
cosmology filled with perfect dust. We will not develop this cosmology
from these properties but merely take it as a given.  In this
cosmology, the universe is modeled as a warped product with leaves
$\Re$ and fibers given by either $\Re^3$, $H^3$ or $S^3$. Thus we have
the metric
\begin{gather}
  ds^2=-d\gt^2 + a^2(\gt)dS^2_K.
\end{gather}
Where $S_K$ is $H^3$, $\Re^3$, or $S^3$ when $K=-1,0,1$, respectively.
The coordinate $\gt$ is time as measured by the isotropic observers.
This is called ``cosmological'' time. The function $a(\gt)$ gives the
scale of the universe, and has the following relationship to $\gt$
depending on $K$
\begin{align}
  K&=1 & a&=\frac{A}2\left(1-\cos u\right) & \gt &=
 \frac{A}2\left(1-\sin u\right)\\
  K&=0 & a&=\left(\frac{9A}4\right)^{1/3}\gt^{2/3}\\
  K&=-1 & a&=\frac{A}2\left(\cosh u-1\right) & \gt &= 
\frac{A}2\left(\sinh u-u\right).
\end{align}
The metric on the fibers is given by
\begin{gather}
  dS^2_K = \frac{dR^2}{1-KR^2}+R^2\left(d\gth^2+\sin^2\gth\, d\gf^2\right).
\end{gather}
If we write 
\begin{gather}
R=  \sin_K(\gc)=\begin{cases}
\sin \gc & \text{ if } \gc = 1\\
 \gc & \text{ if } \gc = 0\\
\sinh \gc & \text{ if } \gc = -1,
\end{cases}
\end{gather}
then we can rewrite the fiber metric as
\begin{gather}
  dS^2_K = d\gc^2 + \sin_K^2\gc\left(d\gth^2+\sin^2\gth\, d\gf^2\right).
\end{gather}
The differences between all of these metrics are small on scales small
compared to the size of the universe, where most of our calculations
will take place.  We also introduce a second time coordinate, $t$, by
the following equation
\begin{gather}
  t=\int\frac{d\gt}{a(\gt)}.
\end{gather}
We can use this to rewrite the metric on the universe as
\begin{gather}
  ds^2(t)= a^2(t)\left(-dt^2+dS^2_K\right).
\end{gather}
For this reason $t$ is called ``conformal time.''

In our discussion of the geometry of a lens system, we will need to
discuss the distances between the observer, lens and source. There are
a number of options available, depending on which equations one wants
to simplify. We will use ``angular diameter distance.'' This distance
is defined by the ratio the physical size of an object to the angular
size of the object seen by the observer.  We denote the distance from
the observe to the lens by $d_L$, the distance from lens to source by
$d_{L,S}$, and the distance from observer to source by $d_{S}$. For
lensing on small scales (nearby galaxies) the universe is almost flat,
and since the bending angle is generally small, we can use $d_S\simeq
d_{L,S}+d_L$.

We will also need the redshift. As the universe expands, the
wavelength of light from distant sources is lengthened. This is
equivalent to clocks appearing to running more slowly. This is
properly associated to an event, but since it changes slowly compared
to the size and duration of a lensing event, we will be assuming it is
locally constant. We can define the redshift of a time, $\gt$, (or an
event) by
\begin{gather}
  z=\frac{a(\gt_O)}{a(\gt)}-1.
\end{gather}

\section{Geometry of Lens System}
We make a number of assumptions about the geometry of a gravitational
lensing system.  First we assume that the lens isn't changing on the
time scale it takes for light to cross the lens. This is valid since
most objects evolve at speeds much less then the speed of light. We
can also assume that the lens is stationary. Any relative motion will
be attributed to the source.
\begin{ass}
  The geometry of the spacetime is assumed to be unchanging on the
  time scale of the lensing event. 
\end{ass}
Furthermore, since almost all sources are ``weak'' we also assume that
the entire system lies in the weak regime, where can approximate the
metric by a time-independent Newtonian potential, $\gf$, given by
\begin{gather}
  \gf(x) = -a^2_L\int_{\Re^3}\frac{\gr(\tilde x)}{\norm{x-\tilde x}}\,d\tilde x.
\end{gather}
Here $a_L$ is the value of $a$ when the light ray is interacting with
the lens.
\begin{ass}
  The spacetime metric of gravitational lens system is given by
  \begin{align}
    ds^2 &= -(1+2\gf)\,d\gt^2+a^2(\gt)(1-2\gf)\,dS_K^2\\
& = a^2(t)\left[-(1+2\gf)\,dt^2+(1-2\gf)\,dS_K^2\right].
  \end{align}
\end{ass}

Furthermore, $\gf$, the Newtonian potential, is much smaller then
unity. We will also assume that while the light ray is interacting
with the lens, $a(t)$ is constant. Furthermore, since on the scale of
this interaction, the universe is almost flat, we will assume $K=0$
during the interaction. Thus near the lens we have the following.
\begin{ass}\Label{ass-flat-potent}
  During the interaction of the light ray with the lens, the metric
  of the lens can be assumed to be 
  \begin{gather}
     ds^2_L  = a^2_L\left[-(1+2\gf)\,dt^2+(1-2\gf)\,dS^2\right].
  \end{gather}
Furthermore, $dt$ can be approximated by $\frac1{a_L}d\gt$.
\end{ass}
We choose dimensionless Euclidean coordinates, $(x_1,x_2,x_3)$,
centered at the lens so that $dS^2 = \gd_{ij}$ and so that the line of
sight to the lens is along the $x_3$ axis. We will also use proper
coordinates $(r_1,r_2,\gz)=a_L\cdot(x_1,x_2,x_3)$, and will use the
coordinates $r=(r_1,r_2)$ on the lens plane. We will use proper
coordinates $s=(s_1,s_2)$ in the source plane. 

\subsection{Mass Densities, Bending Angle and Index of Refraction}
Since we are in the static weak field limit, the Einstein equation
reduces to the time independent Poisson equation for $\gf$. Thus we
get the three dimensional Poisson equation
\begin{gather}
  \gD^{3D} \gf(x) = 4\pi a^2_L \gr(x).
\end{gather}
Where $\gr$ is the density above the background of the lens. The
solution to this equation is
\begin{gather}
  \gf(x) = -a^2_L\int_\Re^3 \frac{\gr\,dx'}{\norm{x-x'}}.
\end{gather}
Since we are assuming that the lens is planar, it is useful to project
the three dimensional potential and density into the lens
plane. Integrating $\gr$ along the line of sight gives us the surface
mass density of the lens $\gs(r)$. This integral is really only from
$-d_{L,S}$ to $d_{L}$, but since we are assuming that $\gr$ is zero
except near the lens we can extend this integral to
$(-\infty,\infty)$. Thus we get that
\begin{gather}
\gs(r) = \int_\Re\gr(r_1,r_1,\gz)\,d\gz.
\end{gather}
Integrating the three dimensional Poisson equation along the $\gz$
axis gives us
\begin{gather}
  4\pi \gs(r) = \int_{\Re} \left(\pp{^2\gf}{r_1^2}+
\pp{^2\gf}{r_2^2}+
\pp{^2\gf}{\gz^2}
\right)\,d\gz
=\gD^{2D}\int_\Re\gf(r_1,r_2,\gz)\,d\gz.
\end{gather}
The $\pp{^2\gf}{\gz^2}$ term integrates to zero since $\gf$ is zero
when $\gz=\pm\infty$.  If we define the surface potential of the lens
by 
\begin{gather}
  \gPs(r) = 2\int_\Re\gf(r_1,r_2,\gz)\,d\gz,\Label{surface-cosmo-pot}
\end{gather}
then $\gPs$ satisfies the two dimensional Poisson equation
\begin{gather}
  \gD^{2D}\gPs(r) = 8\pi\gs(r).
\end{gather}
This is solved by 
\begin{gather}
  \gPs(r) = 4\int_{\Re^2} \gs(r')\ln\norm{\frac{r-r'}{d_0}}\, dr',
\end{gather}
for any constant $d_0$. We will generally choose $d_0=d_L$. 

It will become useful to think of the potential of the lens giving a
refractive index to the spacetime. The index of refraction of a medium
is the reciprocal of the velocity of light in that medium. We want to
calculate the velocity of light in our metric, $ds_L^2$, relative to
the flat metric given by $a^2_L \gd_{ij}$. The velocity of light in
the metric
\begin{gather}
  ds^2  = -A(x)dt^2+B(x)dS^2
\end{gather}
is given by the ratio 
\begin{gather}
  n=\sqrt{\frac{B}{A}}.
\end{gather}
In the metric $ds^2_L$ we get
\begin{gather}
  n=\sqrt{\frac{1-2\gf}{1+2\gf}}\simeq 1-2\gf\Label{index-pot}
\end{gather}
to first order in $\gf$. Since $\gf$ is much smaller then unity, we
can ignore higher order terms. 

Now we can look at the bending angle of the lens. We approximate the
light ray from the source to the observer by a broken null geodesic
with the break at the lens plane. We compress all the bending in the
ray due to the lens into this corner of the light ray. Call the
tangent to the incoming light ray $T_i(r)$ and tangent to the final
ray $T_f(r)$.  Then we define the bending angle by
\begin{gather}
  \hat\ga(r)=T_f(r)-T_i(r).
\end{gather}
Here we have parametrized the vectors by the impact parameter $r$.
Continuing with the standard geometric optics approximation, we
parametrize the spatial path $R(s)=\left(R_1(s), R_2(s), R_3(s)\right)$ of the
light ray by arclength, $s$, in the background metric ($ds^2=a^2_L
\gd_{ij}$). Then light rays are characterized by the equation
\begin{gather}
  \dd{}s\left(n\dd{R}s\right)=\nabla n.\Label{index-eqn}
\end{gather}
Here $\nabla$ is the flat gradient.  Now define the quantities
$T=\dd{R}{s}$ and $K=\dd{T}{s}$ as the tangent and curvature vectors of
the curve $R(s)$.  Plugging these into equation~\Eqref{index-eqn}
gives the equation
\begin{gather}
  (Tn)T + nK =\nabla n.
\end{gather}
Since $K$ is perpendicular to $T$ the transverse gradient is given by
\begin{gather}
  \nabla_\perp n = nK.
\end{gather}
Solving this for $K$ gives us
\begin{gather}
  K = \frac{\nabla_\perp n}n \simeq\left(-2\nabla_\perp\gf\right)
\left(1-2\gf\right)^{-1}\simeq -2\nabla_\perp \gf
\end{gather}
to first order in $\gf$.  Since this angle is small, the light rays
are almost perpendicular to the lens plane, so we can replace
$\nabla_\perp$ with the gradient in the $(r_1,r_2)$ plane, $\nabla_r$,
and we can integrate over $\gz$ to find the total $K$ for the entire
light ray. This total $K$ tells us how far the tangent vector has
turned, hence the bending angle is
\begin{gather}
  \hat\ga(r) = 2\int\nabla_r \gf(r_1,r_2,\gz)\, d\gz. 
\end{gather}
Pushing the integral inside the gradient gives us
\begin{gather}
  \hat\ga(r) = \nabla \gPs(r).\Label{bending-potential}
\end{gather}

\subsection{Fermat's Principle and Time Delays}
\newcommand{\T}{\mathcal T} One could use equation
\Eqref{bending-potential} to try and work out the effect of a lens,
but instead we will follow \Cite{plw} and use the following principle.
\begin{prop}[Fermat's Principle]
  A light ray from a source (an event) to an observer
  (a timelike curve) follows a path that is a stationary value of the
  arrival time functional, $\T$, on paths. 
\end{prop}
Here we are only considering paths $v_r$, that are broken geodesics
from the source $S$, to the observer $O$, parametrized by the impact
parameter $r=(r_1,r_2)$ where the ray crosses the lens plane. For a
given source location $s$ in the source plane, we look at the time
delay function $\T_s(r)$ that gives the time delay for a light ray
that goes from $s$ to $r$, bends at $r$, and then continues on to $O$.
Technically the time delay function and the arrival time function
differ by some reference value. That reference value is the time the
ray would have taken in the absence of the lens. We denote this
unlensed references path by $u_0$. Fermat's principle tells us that
the images of a source at $s$ are given by solutions to the equation
\begin{gather}
  \nabla_r \T_s(r)=0.
\end{gather}

To calculate this we will separate the time delay into two
components. One is the effect of the longer path the light ray
takes. This is called the {\em geometric time delay}, $\T_g$. We will
drop the $s$ for the moment. The other effect is due to time passing
more slowly in a gravitational potential as seen by a distant
observer. This is called the {\em potential time delay}, $\T_p$.  
The travel time for the unlensed ray is given by
\begin{gather}
  \int_{u_0} a_L dl.
\end{gather}
Where $dl$ is given by the metric $dS^2_K$.  The travel time for the
lensed ray is given by
\begin{gather}
  \int_{v_r} a_Ln_L dl.
\end{gather}
Hence the time delay is given by
\begin{gather}
  \T^L(r) = \int_{v_r} a_Ln_L dl-\int_{u_0} a_L dl.
\end{gather}
Here the $L$ attached to $\T$ denotes the fact that these time delays
are being measured at the lens plane. We will have to account for the
redshift, $z_L$, of the lens. Thus we define the geometric and
potential time delays by
\begin{gather}
  \T_p^L = \int_{v_r} a_L(n_L-1) dl\qquad
\T_g^L =  a_L\left(\int_{v_r}dl-\int_{u_0}dl\right).
\end{gather}
The potential time delay is easiest to calculate. Since the bending
angle is small, we can approximate $v_r$ by $u_0$. Using
equations \Eqref{index-pot} and \Eqref{surface-cosmo-pot} we can
calculate the potential time delay as
\begin{gather}
  \T_p^L(r) = -\gPs(r).
\end{gather}
The geometric time delay is more complicated. We will first calculate
the geometric time delay assuming that $K=0$, since we are
calculating many quantities to first order, the geometric time delay
will the same for $K=\pm 1$. For a detailed treatment of $K=\pm 1$,
see \cite{plw}.

First we define the dimensionless lengths $l_S$, $l_L$, and $l_{L,S}$
as the lengths of the spatial projections of $u_0$, and the parts of
$v_r$ between the lens and observer and observer and source
respectively. These lengths are measured in the metric $dS^2_K$. Thus
the geometric time delay measured at the lens is given by
\begin{gather}
  \T_g^L = a_L\left(l_L+l_{L,S}-l_S\right).
\end{gather}
The law of cosines tells us
\begin{gather}
  l_S^2=l_L^2+l_{L,S}^2-2l_Ll_{L,S}\cos(\pi-\hat\ga).
\end{gather}
We can approximate $\cos(\pi-\hat\ga)$ by $-1+\hat\ga^2/2$ since
$\hat\ga$ is small. This gives us
\begin{align}
    l_S^2&\simeq l_L^2+l_{LS}^2+2l_Ll_{L,S}(1-\hat\ga^2/2)\\
&=\left(l_L+l_{L,S}\right)^2-l_Ll_{L,S}\hat\ga^2.
\end{align}
Isolating the term with $\hat\ga$ and factoring gives us
\begin{gather}
  \left(l_L+l_{L,S}-l_S\right)
  \left(l_L+l_{L,S}+l_S\right)\simeq l_Ll_{L,S}\hat\ga^2\\
\T_g^L \simeq\frac{ l_Ll_{L,S}}{l_L+l_{L,S}+l_S}\hat\ga^2
 \simeq\frac{ l_Ll_{L,S}}{2l_S}\hat\ga^2.\Label{tlg-ratio}
\end{gather}
We can replace $l_L$ and $l_{L,S}$ by $a_Ld_L$, and $a_Ld_{L,S}$
respectively. We would like to remove the reference to $\hat\ga$. To
do that we first construct the point $s'$ in the source plane. It is
the location that would produce an image at $r$ in the absence of the
lens. Then using similar triangles we compute
\begin{gather}
  \hat\ga\, d_{L,S} = \norm{s'-s}=\norm{\frac{r}{d_L}-\frac{s}{d_S}}d_S.
\end{gather}
Plugging this into equation \Eqref{tlg-ratio} gives us 
\begin{gather}
  \T_g^L(r)\simeq
  \frac1{a_L}\frac{d_Ld_S}{2d_{L,S}}\norm{\frac{r}{d_L}
    -\frac{s}{d_S}}^2.
\end{gather}
Adding the two parts of the time delay together, and correcting for
the redshift by multiplying by $1+z_L$ gives us
\begin{gather}
  \T_s(r) =
  \left(1+z_L\right)\frac{d_Ld_S}{d_{L,S}}\left(
\frac12\norm{\frac{r}{d_L}-\frac{s}{d_S}}^2
-\frac{d_{L,S}}{d_Ld_S}\gPs(r)
\right).
\end{gather}
Taking the gradient by $r$ and solving for $s$ gives the \emph{lens equation}
\begin{gather}
s=\frac{d_S}{d_L}r-d_{L,S}\hat\ga(r).
\end{gather}
Viewing $s$ as a function of $r$ results in the \emph{lensing map.}
This map goes from the image plane to the source plane. It answers the
question: ``Where would a source have to be create an image at this
location?''  To further simplify this equation, we will
nondimensionalize by introducing the following variables
\begin{gather}
  x=r/d_L \qquad y=s/d_S\\
\gps(x)=\left(\frac{d_{L,S}}{d_Ld_S}\right)\gPs(r)\qquad
\ga(x)=\frac{d_{L,S}}{d_S}\hat\ga(r)
\\
\gk(x)=\frac{\gs(r)}{\gs_c}\text{ where }\gs_c=\frac{d_S}{2\pi
  d_Ld_{L,S}}.
\end{gather}
With these variables, the lensing map becomes
\begin{gather}
y=\gh(x)=x-\ga(x)
\end{gather}
with $\ga=\nabla \gps$. 
\subsection{Magnification}
In addition to changing the location of the image of a source,
gravitational lensing can also change the apparent size of an object.
Since all sources are not truly point sources, we can really consider
how a region, $R$, around a point $s$ in the source plane gets
deformed. In particular, the signed area of the image of $R$ will be
determined by the integral of the determinant of $d\gh^{-1}$. Due to
the Brightness Theorem the apparent surface brightness of an object is
invariant for all observers. For instance if one were twice as far
from the sun, the total light received would be reduced by four, as
would the area of the sun. So the observed surface brightness would
remain constant. Thus the total light received at the observer from an
extended source is scaled the same as the areas. See \cite{MTW} for
more information.

However, while many sources aren't point sources they are point-like.
Thus the only effect of the magnification is to increase (or decrease)
the brightness of the source. Thus the magnification of a point source
at the location $y$ is given by
\begin{gather}
  \gm_y(x) =\frac{1}{\abs{\det d\gh(x)}}.
\end{gather}
Locations, $x$, in the lens plane where this is infinite are called
{\em critical points}. The corresponding locations $\gh(x)=y$ in
source plane are called {\em caustics}. Sources on different sides of
caustics typically have two fewer or less images then each other. It
is often useful to consider a source that is moving in the source
plane. As this source crosses the caustic, two of its images will
increase in brightness and merge, then disappear, or the reverse
depending on the direction in which the source crosses the caustic.
The changing magnification of a source as it moves in the source plane
is a useful observable called the {\em light curve.} For these curves,
we add up the magnification of all the images, since sometimes the
images are too close to be resolved.

\section{Lensing Map for Isolated Negative Point Mass Singularities}
The general framework we have established for gravitational lensing
requires a potential function to plug into the metric in
assumption~\Ref{ass-flat-potent} and following formulas. We will be
studying singularities that agree to first order with negative
Schwarzschild solutions.  This is summarized in the following
assumption:
\begin{ass}
  The two dimensional surface mass density for a negative point mass
  singularity is given by
\begin{gather}
\gr(r) = M\gd(r).
\end{gather}
\end{ass}
To first order, this assumption gives the correct behavior of a
negative mass Schwarzschild solution in the weak field regime. We will
restrict ourselves to cases that agree with this case to first
order.\footnote{For a more detailed assumption of the positive mass
  analogue of this assumption see the discussion of point masses on
  p.\ 101 in \Cite{plw}.}

We nondimensionalize our potential by defining
\begin{gather}
  m = \frac{M}{\gp d^2_L\gs_c}.
\end{gather}
Which gives us the dimensionless surface mass density
\begin{gather}
  \gk(x) = \gp m \gd(x).
\end{gather}
This gives us the dimensionless surface potential
\begin{gather}
  \gps(x)=m\ln\left(\norm{x}\right),
\end{gather}
and the dimensionless bending angle
\begin{gather}
  \ga(x) = m\frac{x}{\norm{x}^2}.
\end{gather}
Thus the dimensionless lensing map is given by
\begin{gather}
  y=\gh(x) = x\left(1-\frac{m}{\norm{x}^2}\right).
\end{gather}
In this case, we can solve this equation exactly to find the images of
a source at $y$. 
\begin{gather}
  x_\pm = \frac12\left(\norm{y}\pm \sqrt{\norm{y}^2+4m}\right)\hat y.
\end{gather}
Here $\hat y$ is the unit vector in the direction of $y$ from the
origin. As long as $\norm{y}>2\sqrt{-m}$, we get two images for each
source.  Both images are on the same side of the lens. If we had
$m>0$, then our images would be on opposite sides of the lens. 
The derivative of the lensing map is given by
\begin{gather}
\begin{pmatrix}
1-\frac{m}{\norm{x}^2}+\frac{2x_1^2m}{\norm{x}^4}&
\frac{2x_1x_2m}{\norm{x}^4}\\
\frac{2x_1x_2m}{\norm{x}^4}&
1-\frac{m}{\norm{x}^2}+\frac{2x_2^2m}{\norm{x}^4}
\end{pmatrix}.
\end{gather}
Hence the magnification is given by
\begin{gather}
\mu(x)=\frac1{1-\frac{m^2}{\norm{x}^4}}=\frac{\norm{x}^4}{\norm{x}^4-m^2}.
\end{gather}
For $x_-$ this number is negative, so to get the total magnification
we take the difference between signed magnifications\footnote{See
  Appendix \Ref{append-calc-mag} for the calculation.}
\begin{gather}
  \gm_{\text{tot}}(y) = \gm_y(x_+)-\gm_y(x_-) = 
\frac{\norm{y}^2+2m^2}{\norm{y}\sqrt{\norm{y}^2+4m}}.
\end{gather}

\section{Light Curves for Isolated Negative Point Mass Singularities}
As we calculated, the lensing map is given by 
\begin{gather}
  y=\gh(x) = x\left(1-\frac{m}{\norm{x}^2}\right).
\end{gather}
This has inverse 
\begin{gather}
  x_\pm = \frac12\left(\norm{y}\pm \sqrt{\norm{y}^2+4m}\right)\hat y.
\end{gather}
The inverse map tells us where an image will appear for a source
located at $y$ in the source plane. Since the radical is imaginary for
$\norm{y}<2\sqrt{-m}$, these sources aren't visible at all. For
sources outside this disk, we get two images: $x_\pm$. The reversed
image $x_-$ is closer to the center of the lens. For large $\norm{y}$,
the $x_-$ image is closer and closer to the center of the lens, while
the $x_+$ image is closer and closer to its unlensed location.  As $y$
gets closer and closer to the caustic $\norm{y}=2\sqrt{-m}$, the two
images come together at $x=\sqrt{-m}$. Looking at the individual
magnifications, we see that for large $\norm{y}$, the positive image
has magnification about 1, while the negative image has magnification
about $-m^2/\norm{y}^4$. As $y$ gets closer and closer to the caustic,
the magnification of each image increases, and is formally infinite
when $y=2\sqrt{-m}$ and $x\pm= \sqrt{-m}$. This curve is the critical
curve. When $y$ is inside the caustic, it creates no image.

In cases where the source, lens and observer are moving relative to
each other, the total magnification changes. We look at the total
magnification, since the individual images are unresolvable, and the
only effect of the lensing is the magnification. By the symmetry of
the source, these paths are characterized by impact parameter, the
distance of closest approach to $y=0$. If this distance is $d$, and we
assume that the source is moving at constant unit speed, the total
magnification as a function of time is
\begin{gather}
  \gm_{m,d}(t) =
  \frac{d^2+t^2+2m}{\sqrt{d^2+t^2}\sqrt{d^2+t^2+4m}}.\Label{light-curve-formula}
\end{gather}
Figure~\Ref{micro-path} shows several possible paths for a moving
source.
\begin{figure}
\begin{center}
\includegraphics[width=2in]{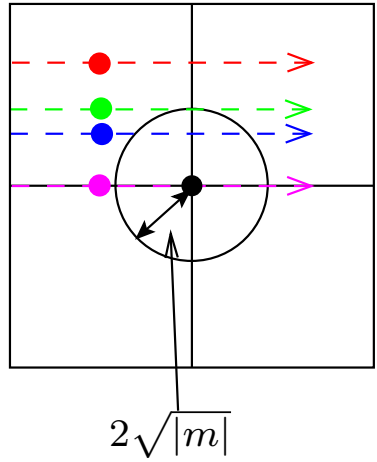}
\caption{Source Paths for Negative Point Mass Microlensing}
\Label{micro-path}
\end{center}
\end{figure}
The sources in these curves have impact parameters varying from zero
to twice $2\sqrt{-m}$.  Figure~\Ref{micro-mag} shows the corresponding light curves. 
\begin{figure}
\begin{center}
\includegraphics[width=5in]{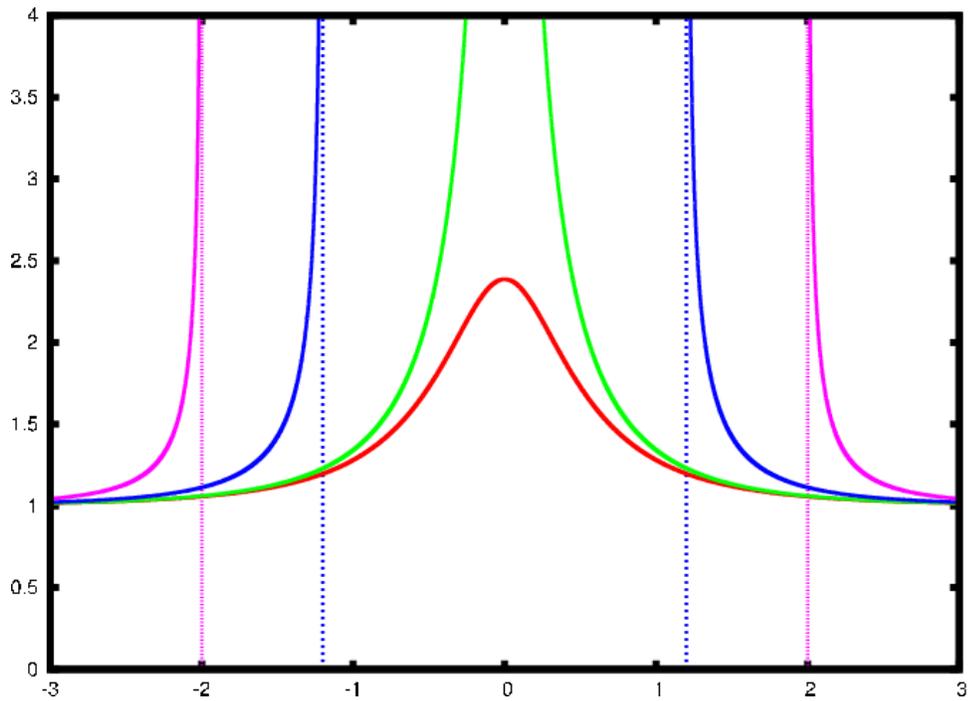}
\caption{Light Curves for Negative Point Mass Microlensing}
\Label{micro-mag}
\end{center}
\end{figure}
The light
curves for the sources passing inside of $2\sqrt{-m}$ are distinctive,
but those for the source passing outside or just along the caustic are
not. The light curve formula, equation~\Eqref{light-curve-formula}, is
the same as that for a positive point mass. Thus if we define 
\begin{gather}
\tilde m= -m\qquad\qquad\tilde d = \sqrt{d^2+4m},  
\end{gather}
then the light curve for a source passing within $\tilde d$ of the
line of sight of a mass of $\tilde m$ is the same as that for a source
passing within $d$ of the line of sight of a mass of $m$ 
\begin{align}
  \gm_{\tilde m,\tilde d}(t) &= \frac{\tilde d^2+t^2+2\tilde
    m}{\sqrt{\tilde d^2+t^2}\sqrt{\tilde d^2+t^2+4\tilde m}}\\
&=\frac{d^2+4m+t^2-2m}{\sqrt{d^2+4m+t^2}\sqrt{d^2+4m+t^2-4m}}\\
&=\gm_{ m, d}(t).
\end{align}
Furthermore, we will see that if we introduce continuous matter we can
reproduce the entire lensing map with a positive mass singularity.

\section{Complex Formulation}
Before incorporating additional features into our lens, it is helpful
to reframe the structure of the lensing map as a map $\Ce\to\Ce$,
rather than a map $\Re^2\to\Re^2$. First we will consider $x$ as a
complex number $x_1+ix_2$, and likewise $y$ and $\gh$. 
The lens equation becomes
\begin{gather}
  \gh=\gh_1 + i\gh_2 = \left(x_1-\pp{\gps}{x_1}\right)+
  i\left(x_2-\pp{\gps}{x_2}\right)
\end{gather}
Taking complex derivatives of $\gh$ we get
\begin{multline}
  \pp{\gh}{z} = \frac12\left(\pp\gh{x_1}-i\pp\gh{x_2}\right)=\\
  \frac12\left(\pp{\gh_1}{x_1}+i\pp{\gh_2}{x_1}-i\pp{\gh_1}{x_2}
    +\pp{\gh_2}{x_2}\right)
  =  \frac12\left(\pp{\gh_1}{x_1}+\pp{\gh_2}{x_2}\right),
\end{multline}
which is real. Here we used the equality of the mixed partials of
$\gps$ to cancel the imaginary terms.  Differentiating with respect to
$\bar z$ gives
\begin{gather}
  \pp{\gh}{\bar{z}} = \frac12\left(\pp\gh{x_1}+i\pp\gh{x_2}\right)=
  \frac12\left(\pp{\gh_1}{x_1}-\pp{\gh_2}{x_2}+i\left[\pp{\gh_2}{x_1}+
    \pp{\gh_1}{x_2}\right]\right).
\end{gather}
Here no such cancellation occurs. 
We can also rewrite $J=\det(d\gh)$ as 
\begin{gather}\begin{aligned}
  J &= \pp{\gh_1}{x_1}\pp{\gh_2}{x_2}- \pp{\gh_1}{x_2}\pp{\gh_2}{x_1}\\
&=\frac14 \pp{\gh_1}{x_1}^2+\frac12\pp{\gh_1}{x_1}\pp{\gh_2}{x_2}+
\frac14\pp{\gh_2}{x_2}^2-
\left(\frac14  \pp{\gh_1}{x_1}^2-
\frac12\pp{\gh_1}{x_1}\pp{\gh_2}{x_2}
+\frac14\pp{\gh_2}{x_2}^2
\right)\\
&\qquad-\left(\frac14  \pp{\gh_2}{x_1}^2+
\frac12\pp{\gh_2}{x_1}\pp{\gh_1}{x_2}
+\frac14\pp{\gh_1}{x_2}^2
\right)\\
&=\abs{\pp{\gh}{z}}^2-\abs{\pp{\gh}{\bar{z}}}^2.
\end{aligned}
\end{gather}
Here we extensively used the fact that
$\pp{\gh_2}{x_1}=\pp{\gh_1}{x_2}$. Our critical points are located
where $J=0$. As we noted above $\pp{\gh}{z}$ is real so the
solutions to $J=0$ look like
\begin{gather}
  \pp{\gh}{\bar z } = \abs{\pp{\gh}{z}}e^{i\gf},\Label{crit-param}
\end{gather}
for some angle $\gf$. Thus our critical curves will be curves
parametrized by $\gf$.

The caustics will be the images of these critical curves under
$\gh$. Any points where the caustics aren't smooth are characterized by
\begin{gather}
  J(x) = 0 \qquad \grad_Z(\gh) = 0. \Label{cusp-complex}
\end{gather}
Here $Z = -\pp{J}{x_2}+i\pp{J}{x_1}=2i\pp{J}{\bar z}$, and
\begin{gather}
  \grad_Z = Z\pp{}z+Z\pp{}{\bar z}.
\end{gather}
To find the cusps on the caustics we just find the appropriate phase
$\gf$ to solve \Eqref{cusp-complex}.

\section{Lensing by Negative Point Mass Singularities with 
Continuous Matter and Shear}
Most lensing events on the scale of individual stars take place within
a host galaxy. In these cases, the star itself isn't the only source
of distortion. Two other non-local factors are also important.

Continuous matter is the first.  The presence of evenly dispersed
matter in the area of the lens can produce convergence or divergence.
This enters into the dimensionless potential via a term like
\begin{gather}
  \gps_{\text{cm}}(x) = \frac\gk2\norm{x}^2. 
\end{gather}
This gives us a lensing map of 
\begin{gather}
 \gh(x)= y = (1-\gk)x.
\end{gather}
This is clearly compatible with the complex formulation.  The
dimensionless mass density of $\gk$ corresponds to a surface mass
density of $\gs_c\gp\gk$.

The other factor is shear from infinity. The presence of a large mass
nearby, such as a nearby galaxy, or an asymmetric distribution such as
the disk of the host galaxy, can introduce a potential of the form
\begin{gather}
  \gps_{\text{sh}}(x) =
  -\frac\gg2\left[\left(x_1^2-x_2^2\right)\cos2\gth+2x_1x_2\sin 2\gth\right].
\end{gather}
The parameter $\gg$ defines the magnitude of the shear, and $2\gth$
determines the preferred direction of the asymmetric mass distribution
or the direction to the large mass. By the symmetry of our lens, we
can assume that $\gth =0$. This gives us a lensing map of
\begin{gather}
   y=
  \begin{bmatrix}
    1+\gg & 0 \\ 0 & 1-\gg
  \end{bmatrix}x.
\end{gather}
In complex form this is
\begin{gather}
  y = z+\gg\bar z.
\end{gather}

Now we will incorporate all of these features into a single lens. The
lensing map of a negative point mass lens is given by
\begin{gather}
  y =  x\left(1-\frac{m}{\norm{x}^2}\right).
\end{gather}
In complex form this is
\begin{gather}
  y = z -\frac{m}{\bar z}.
\end{gather}
Adding these potentials to that of our singularity gives us the
combined lens equation
\begin{gather}
  \gh(x) = \begin{bmatrix}1-\gk+\gg& 0 \\ 0 & 1-\gk-\gg\end{bmatrix}x - \frac{m}{\norm{x}^2}x,
\end{gather}
or
\begin{gather}
  \gh = (1-\gk)z + \gg\bar z -\frac{m}{\bar z}.
\end{gather}
Our Jacobian is
\begin{gather}
  J= \left(\gg +\frac{m}{\bar z^2}\right)^2- \left(1-\gk\right)^2.
\end{gather}
So to look for critical points we use equation \Eqref{crit-param} with
our $\gh$ to give us
\begin{gather}
  \gg +\frac{m}{\bar z^2}= \abs{1-\gk}e^{i\gf}.
\end{gather}
If $\gk=1$, then we are looking for points that solve 
\begin{gather}
  \frac{m}{\bar z^2} = \gg.
\end{gather}
Which are just two points along the $x$ axis. So we can assume that
$\gk\neq 1$. It simplifies calculation to remove $\gk$ from the
calculation by defining
\begin{gather}
  \gg_* = \frac{\gg}{\abs{1-\gk}}\qquad
  m_* = \frac{m}{\abs{1-\gk}}
\qquad \ge_\gk = \operatorname{sgn}(1-\gk).
\end{gather}
With those substitutions we now solve
\begin{gather}
  \gg_* + \frac{m_*}{\bar z^2 } = e^{i\gf}.
\end{gather}
This has solutions
\begin{gather}
  z_{\pm}(\gf) = \pm\sqrt{\frac{m_*}{e^{-i\gf}-\gg_*}}.
\end{gather}
At this point the only effect that the negative sign on $m$ has had is
to rotate our critical curves by $\pi/2$. 
\begin{gather}
  z_{\pm}(\gf) =  \pm i\sqrt{\frac{\abs{m_*}}{e^{-i\gf}-\gg_*}}.
\end{gather}
Note that these curves are independent of $\ge_\gk$.  For $\gg_*\ll1$,
the critical curve is given by an oval with long axis in the $x_2$
direction. As $\gg_*$ grows to $1$, the critical curve develops a
waist.  For $\gg_*>1$, we have two critical curves, small loops on the
$x_1$ axis.  As $\gg_*$ grows they shrink to the two points for
$\gk=1$.  For $\gg_*$ very close to $1$, both varieties of critical
curves grow to infinity. When $\gg_*=1$ our critical curves degenerate
into two curves asymptotic to the lines $x_2=\pm x_1$.  As $\gg_*$
passes $1$, the ends of the curve open up, pass through infinity, and
reseal re-paired.  See Figure~\Ref{crit-fig} for the shapes of the
critical curves for various values of shear and convergence. 
\begin{figure}
\begin{center}
  \begin{tabular}{ccc}
\includegraphics{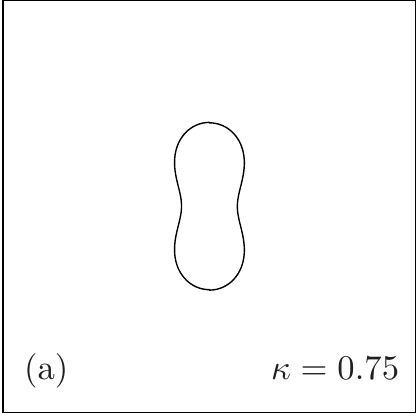}&
\includegraphics{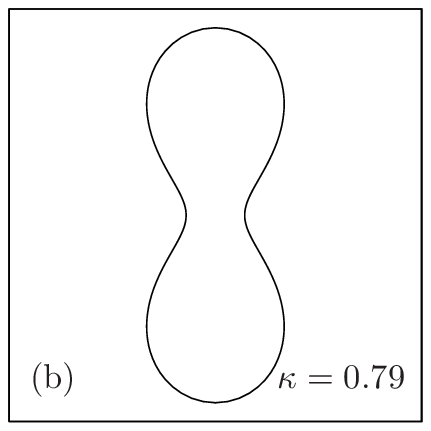}&
\includegraphics{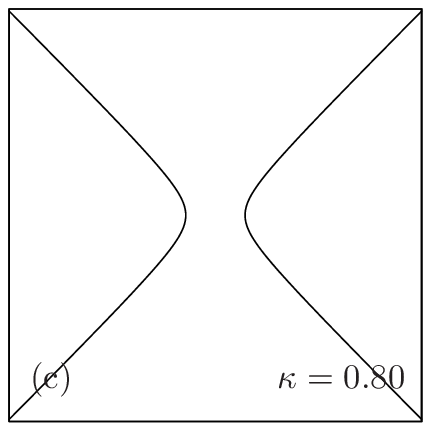}\\
\includegraphics{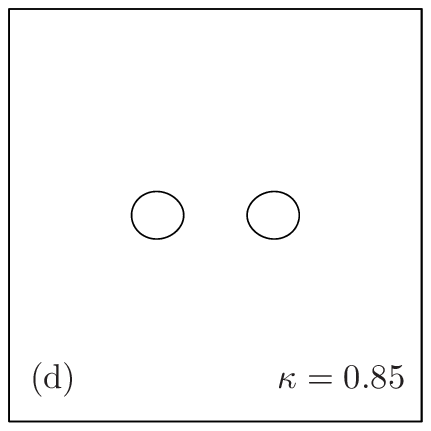}&
\includegraphics{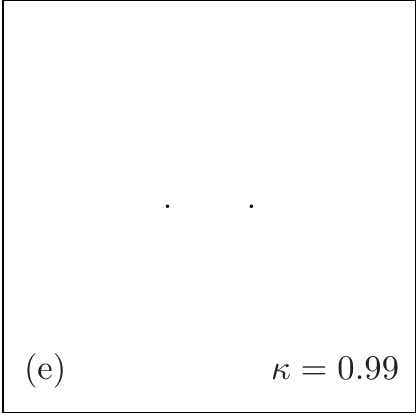}&
\includegraphics{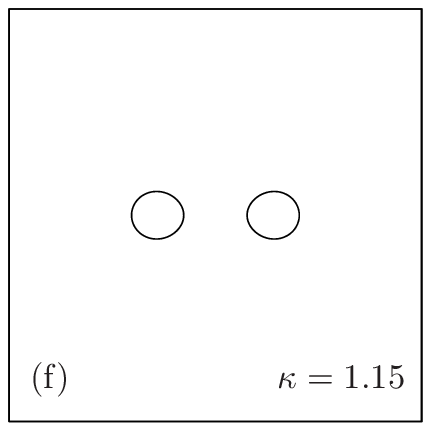}\\
\includegraphics{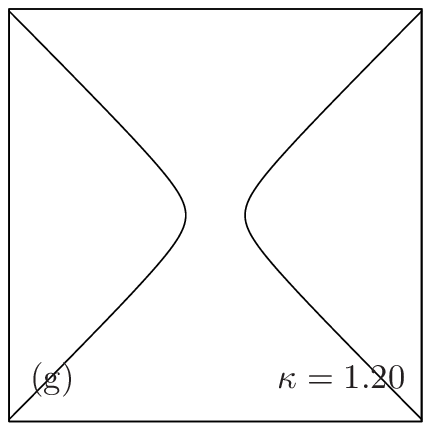}&
\includegraphics{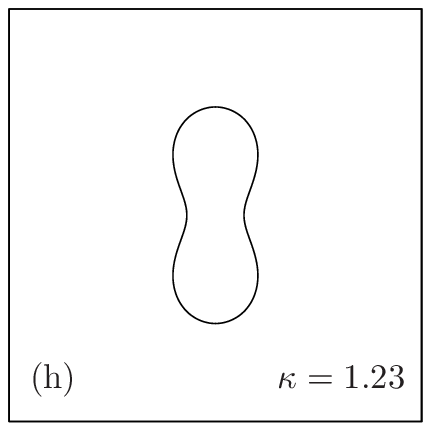}&
\includegraphics{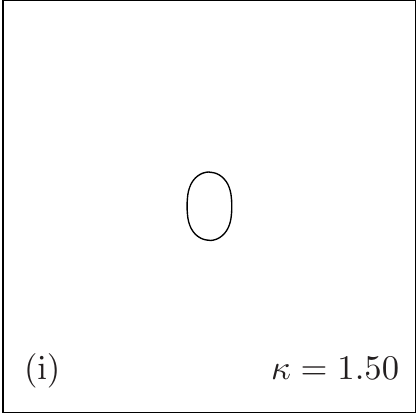}
\end{tabular}
  \caption{Critical Curves for Negative Point Mass with Continuous
    Matter and Shear}\Label{crit-fig}
\end{center}

These pictures are shown in order of increasing $\gk$, with constant
$\gg=0.2$. The first four correspond to $\ge_\gk=1$ and increasing
$\gg_*$. At $\gk=1$ we have two points as the critical set. The last
four pictures correspond to $\ge_\gk=-1$ and decreasing $\gg_*$. Since
the curves are independent of $\ge_\gk$, they are the same when one
replaces $\gk$ by $2-\gk$.
\end{figure}
These curves are the same as in the positive mass case, rotated a
quarter turn.

To find the cusps we plug our lensing map into \Eqref{cusp-complex}.
In our case we have
\begin{gather}
 Z= -4i\left(\gg+\frac{m}{z^2}\right)\frac{m}{\bar z^3}
\end{gather}
and equation \Eqref{cusp-complex} is 
\begin{gather}
  0 = (1-\gk)Z + \left(\gg+\frac{m}{\bar z^2}\right)\bar Z. 
\end{gather}
Removing the $1-\gk$ as before, we get
\begin{gather}
   0 = Z_*  +\ge_\gk\left(\gg_*+\frac{m_*}{\bar z^2}\right)\bar Z_*. 
\end{gather}
Where $Z_*$ has the same formula as $Z$ replacing $\gg$ and $m$ with
$\gg_*$ and $m_*$.  In Appendix~\Ref{append-calc-cusp}
we find the possible roots of this equation to be
\begin{gather}
  \gf_1 = 0,\quad \gf_2=\pi, \quad \gf_{3,4} =\acos\left( \frac{3\pm
      \sqrt{4\gg_*^2-3}}{4\gg_*} \right),\quad \gf_5 =
  2\pi-\gf_3,\quad \gf_6 = 2\pi-\gf_4.
\end{gather}
For a given value of $\gg_*$ and $m_*$, each of these is only a
solution for either $\ge_\gk=1$ or $\ge_\gk=-1$. The numbers of cusps for
various values of $\gg_*$ and $\ge_\gk$ are given in
Table~\Ref{cusp-table}.
\begin{table}
  \centering
  \begin{tabular}{|c|cc|cc|}
    \hline
    &\multicolumn{2}{c|}{$\ge_\gk=1$}&
    \multicolumn{2}{c|}{$\ge_\gk=-1$}\\
    Shear &$\gf_i$&$N_\text{cusps}$&$\gf_i$&$N_\text{cusps}$
    \\\hline\hline 
    $0\leq\gg_*^2<3/4$ & $\emptyset$ & 0 & $\gf_1$, $\gf_2$ & 4 \\\hline
    $3/4\leq\gg_*^2<1$ & $\gf_3$, $\gf_4$, $\gf_5$, $\gf_6$ & 8
    & $\gf_1$, $\gf_2$ & 4 \\\hline
    $1<\gg_*^2$ & $\gf_1$, $\gf_4$, $\gf_6$ & 6
    & $\gf_2$, $\gf_3$, $\gf_5$ & 6 \\\hline
 \end{tabular}
  \caption{Numbers of cusps on caustics for various shear and continuous matter values.}\Label{cusp-table}
\end{table}
 When the number of cusps changes we expect to
have higher order caustics. Looking at Figure~\Ref{caustic-fig},
\begin{figure}
 \begin{center}
  \begin{tabular}{ccc}
\includegraphics{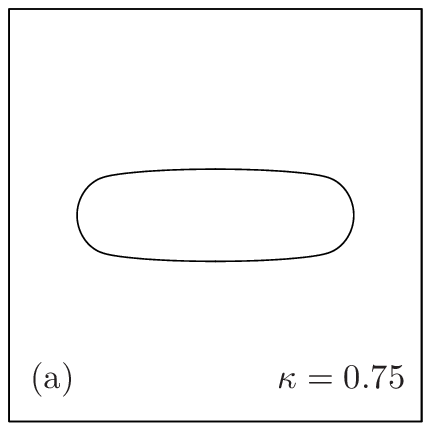}&
\includegraphics{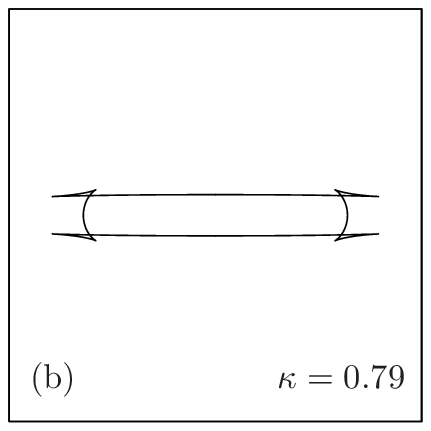}&
\includegraphics{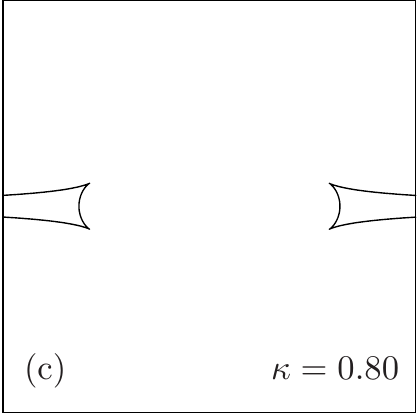}\\
\includegraphics{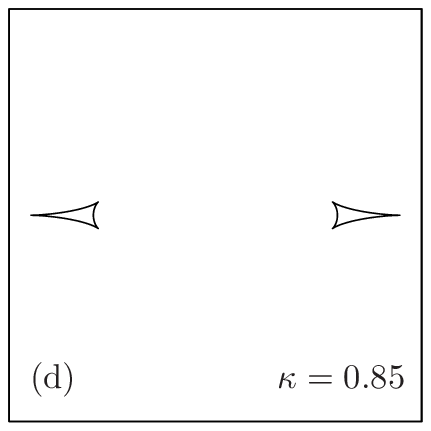}&
\includegraphics{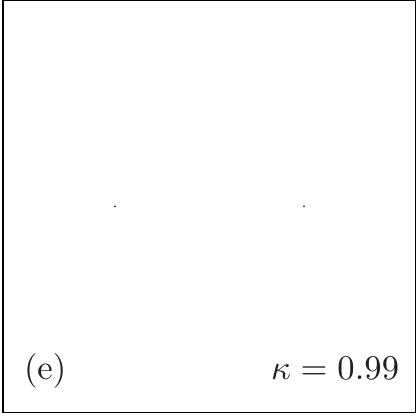}&
\includegraphics{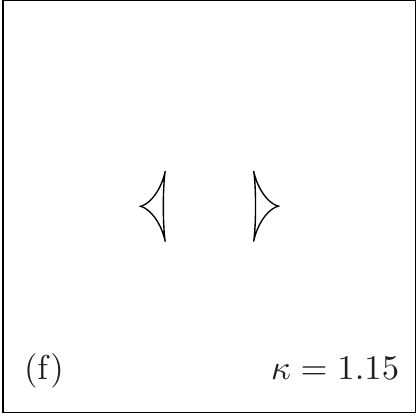}\\
\includegraphics{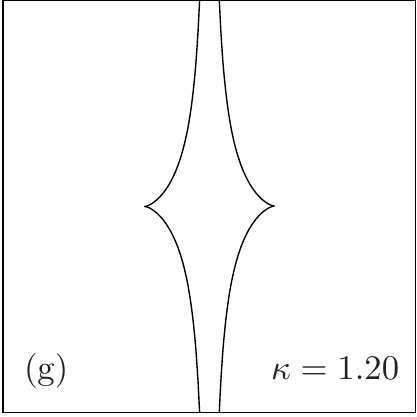}&
\includegraphics{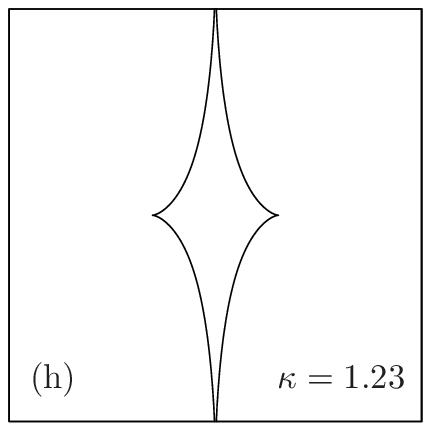}&
\includegraphics{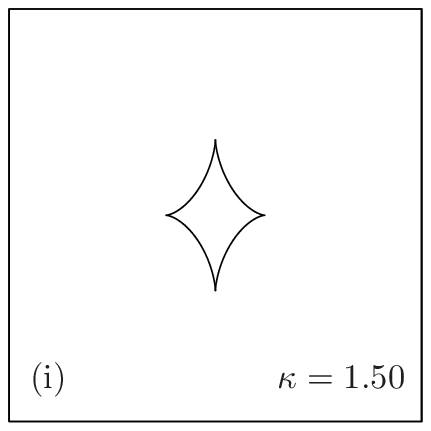}
\end{tabular}
  \caption{Caustics for Negative Point Mass with Continuous
    Matter and Shear}\Label{caustic-fig}
\end{center}

These pictures are shown in order of increasing $\gk$ with constant $\gg=0.2$. The first four
correspond to $\ge_\gk=1$ and increasing $\gg_*$. At $\gk=1$ we have
two points as the caustic set. The last four pictures correspond to
$\ge_\gk=-1$ and decreasing $\gg_*$. 
 \end{figure}
we see that we have four swallow tails between panels $(a)$ and $(b)$,
we also have two elliptic umbilics between panels $(d)$ and $(f)$.
The changes in the caustic structure between panels $(b)$ \& $(d)$ and
$(f)$ \& $(h)$ occur at infinity. 

In terms of the number of images, we always get two images outside of
the caustics. Inside the caustic in panel $(a)$ we have no images,
just like in the case without shear or convergence. In panel $(b)$ we
have four images inside the swallow tails, and none inside the
rectangular region. In the rest we have four images inside the
caustics. 

Both the critical curves and the caustics are the same as what one
would get if one took the positive mass case, and switched $\gk$ with
$2-\gk$, and rotated all the images by a quarter turn.

\chapter{Inverse Mean Curvature Flow}\Label{IMCF}
This chapter lays out what we need of the weak inverse mean curvature
flow as developed in \cite{IMCF}. We will follow their exposition
closely, omitting many of the technical details.

\section{Classical Formulation}
Let $N$ be a the smooth boundary of a region in the smooth Riemannian
manifold $M$. A classical solution of the inverse mean curvature flow
is a smooth family $x:N\times[0,T]\to M$ of hypersufaces $N_t=x(N,t)$
satisfying the evolution equation
\begin{gather}
  \pp{x}{t} = \frac{\gv}{H}.\Label{IMCF-classic}
\end{gather}
Here $\gv$ is the outward pointing normal to $N_t$ and $H$ is the mean
curvature of $N$, which must be positive. We use this flow to explore
the geometry near a singularity. For instance we will use the fact
that under this flow the Hawking mass is non-decreasing.  This result
can be readily shown if we assume that the flow doesn't have any
discontinuities or singularities. However, easy counterexamples
illustrate that this is overly optimistic. The simplest counterexample
is given by a thin torus in $\Re^3$. Such a torus has positive mean
curvature approximately that of a cylinder of the small radius. The
inverse mean curvature flow will tend to increase the small radius of
this torus. However, if it continued without singularities or jumps,
the hole in the torus would eventually shrink to the point where the
area of the torus inside the hole has zero mean curvature and hence
the flow couldn't be continued.  Thus the classical flow is
insufficient for our needs.

To remedy this problem, we will follow \Cite{IMCF} and recast the flow
first in a level set formulation and then we will move to a weak
solution. This will allow the flow to jump to avoid situations where
the curvature of the surface would drop to zero. In the previous
example, the flow would close the interior of the torus as soon as it
is favorable in terms of a certain energy functional. Even with these
jumps the Hawking mass of our surface is still non-decreasing.

\section{Weak Formulation}
First we establish some notation. Let $(M,g)$ be the ambient manifold.
Let $h$ be the induced metric on $N$. Let $A_{ij}=\ip{\nabla_{e_i}
  \gv}{e_j}$ be the second fundamental form of $N$. Then $H$ is the
trace of $A$ with $h$, and $\vec H=-\gv H$ is the mean curvature
vector. Let $E$ be the open region bounded by $N$.

The first step toward the weak formulation is a level set
formulation. We assume that the flow is given by the level sets of a
function $u:M\to\Re$. This $u$ is related to our previous data by
\begin{gather}
  E_t:=\left\{x\mid u(x)<t\right\}, \qquad N_t:=\p E_t.
\end{gather}
We will also need the following sets
\begin{gather}
  E_t^+:=\text{int}\left\{x\mid u(x)\geq t\right\}, \qquad
  N_t^+:=\p E^+_t
\end{gather}
Where $\grad u\neq0$, $E_t=E_t^+$ and $N_t=N_t^+$. 

Anywhere that $u$ is smooth and $\nabla u\neq0$, then we have a
foliation by smooth surfaces $N_t$, with normal vector $\gv=\grad
u/\abs{\grad u}$. The mean curvature of these surfaces is given by
$\Div_N(\gv)$ and the flow velocity is given by $1/\abs{\grad{u}}$,
so equation \Eqref{IMCF-classic} becomes
\begin{gather}
  \Div_M\left(\frac{\grad u}{\abs{\grad u}}\right)= \abs{\grad
    u}.\Label{IMCF-weak}
\end{gather}
This equation is degenerate elliptic. To remedy this we introduce the
functional $J_u^K(v)$:
\begin{gather}
  J_u(v)=J_u^K(v)=\int_K \abs{\grad v}+v\abs{\grad u}.
\end{gather}
Where $K$ is a compact set in $M$.  If we take the Euler-Lagrange
equation of this functional, and replace $v$ with $u$, we get back
equation \Eqref{IMCF-weak}. For each $u$ we get a different $J_u$.
Hence what we want is a function $u$ which minimizes its own $J_u$.
\begin{defn}
  Let $u$ be a locally Lipschitz function on the open set $\gO$ in
  $M$. Then $u$ is a \emph{weak (sub-, super-) solution} of
  equation \Eqref{IMCF-weak} on $\gO$ exactly when 
  \begin{gather}
    J_u^K(u)\leq J_u^K(v)
  \end{gather}
  for all locally Lipschitz functions $v$ ($\leq u,\geq u$) which only
  differ from $u$ inside a compact set $K$ contained in $\gO$. 
\end{defn}
It is worth noting that $J_u(\min(v,w))+J_u(\max(v,w)) = J_u(v)+J_u(w)$.
To see this, construct $K_v= \left\{x\in K\mid v(x)<w(x)\right\}$, and
divide the integrals on the left into integrals over $K_v$ and
$K\setminus K_v$.  Then regroup then and recombine to get the right
hand side.  This tells us that if $u$ is a both a weak supersolution
and subsolution, then it is a weak solution.

We will also need a related functional of the level sets. 
\begin{defn}If $F$ is a
set of locally finite perimeter, and $\p^*F$ is its reduced boundary,
then we define 
\begin{gather}
  J_u(F) = J_u^K(F)=\abs{\p^*F\cap K}-\int_{F\cap K}\abs{\grad u}.
\end{gather}
For any locally Lipschitz function $u$ and compact $K$ contained in
$A$. We say that $E$ \emph{minimizes $J_u$ in $A$ (on the inside,
  outside)} if
\begin{gather}
J_u^K(E)\leq J_u(F)
\end{gather}
for all $F$ that differs from $E$ in some compact $K$ in $A$ (with
$E\subseteq F$, $E\supseteq F$.)  A similar argument tells us that if
$E$ minimizes $J_u$ exactly when it minimizes $J_u$ on the inside and
outside. 
\end{defn}
These two formulations are equivalent.
\begin{lemma}
  Let $u$ be a locally Lipschitz function in the open set $\gO$. Then
  $u$ is a weak (sub-, super-) solution of equation \Eqref{IMCF-weak}
  exactly when for each $t$, $E_t=\left\{u<t\right\}$ minimizes $J_u$
  in $\gO$ (on the inside, outside).
  \begin{proof}
    Lemma 1.1 in \cite{IMCF}.
  \end{proof}
\end{lemma}
We now define the initial value problem.  Let $E_0$ be an open set
with $C^1$ boundary. We say that $u\in C_\text{loc}^{0,1}$ and the
associated $E_t$ for $t>0$ is a \emph{weak solution of
  \Eqref{IMCF-weak} with initial condition $E_0$} if either
 \begin{gather}
  \begin{gathered}
   E_0=\{u<0\}\text{ and }u\text{ minimizes }J_u\text{ on }M\setminus E_0 \\
\text{or}\\
   E_t=\{u<t\}\text{  minimizes }J_u\text{ in }M\setminus E_0\text{ for }t>0.
   \end{gathered}\Label{IMCF-IVP}
 \end{gather}

These two conditions are equivalent by Lemma 1.2 in \cite{IMCF}. 
Showing the regularity of $N_t$ and $N_t^+$ is nontrivial, but we
won't reproduce it here. 
\begin{thm}
  Let $n<8$. Let $U$ be an open set in a domain $\gO$. Let $f$ be a
  bounded measurable function on $\gO$. Consider the functional
  \begin{gather}
    \abs{\p F}+\int_Ff
  \end{gather}
  on sets containing $U$ and compactly contained in $\gO$. Suppose $E$
  minimizes this functional. 
  \begin{enumerate}
  \item If $\p U$ is $C^1$, then $\p E$ is a $C^1$ submanifold of $\gO$.
  \item If $\p U$ is $C^{1,\ga}$, $0<\ga\leq 1/2$, then $\p E$ is a
    $C^{1,\ga}$ submanifold of $\gO$. The $C^{1,\ga}$ estimates depend
    only on the distance to $\p\gO$, $\operatorname{ess}\sup\abs{f}$,
    $C^{1,\ga}$ bound for $\p U$, and $C^1$ bounds on the metric in
    $\gO$.
  \item If $\p U$ is $C^2$ and $f=0$, then $\p E$ is $C^{1,1}$, and
    $C^\infty$ where it doesn't touch $U$. 
  \end{enumerate}
\end{thm}
Our initial value formulation falls into this category of problem. So
our $N_t$'s and $N_t^+$'s are $C^{1,\ga}$. Furthermore 
\begin{gather}
  \lim_{s\to t^-} N_s=N_t\qquad \lim_{s\to t^+}N_s=N_t^+
\end{gather}
in local $C^{1,\gb}$ convergence $0<\gb\leq\ga$. 

The locations where $N_t\neq N_t^+$ correspond to the jumps discussed
earlier. To examine these areas we need to introduce minimizing hulls.
\begin{defn}
  Let $\gO$ be an open set. We call $E$ a \emph{minimizing hull} if
  \begin{gather}
    \abs{\p^* E \cap K} \leq \abs{\p^* F\cap K}
  \end{gather}
  for any $F$ containing $E$ with $F\setminus E$ in $K$ a compact set
  in $\gO$. We say $E$ \emph{strictly minimizes} if equality implies
  that $E$ and $F$ agree in $\gO$ up to measure zero. 
\end{defn}
The intersection of (strictly) minimizing hulls is a (strictly)
minimizing hull. So, given a set $E$ we can intersect all of the
strictly minimizing hulls which contain $E$. This gives, up to measure
zero, a unique set $E'$ that we will call \emph{the strictly
  minimizing hull of $E$}. Since $E'$ is strictly minimizing,
$E''=E'$.

Solutions to the initial value problem given by equation 
\Eqref{IMCF-IVP}, have level sets that are minimizing hulls as follows.
\begin{lemma}
  Suppose that $u$ is a solution to \Eqref{IMCF-IVP}. Assume that $M$
  has no compact components. Then:
  \begin{itemize}
  \item For $t>0$, $E_t$ is a minimizing hull in $M$.
  \item For $t>0$, $E_t^+$ is a strictly minimizing hull in $M$.
  \item For $t>0$, $E_t' = E_t^+$ if $E_t^+$ is precompact. 
  \item For $t>0$, $\abs{\p E_t} = \abs{\p E_t^+}$, provided that
    $E_t^+$ is precompact. 
  \item Exactly when $E_0$ is a minimizing hull $\abs{\p E_0} = \abs{\p E_0^+}$
  \end{itemize}
\end{lemma}
This is ``Minimizing Hull Property (1.4)'' in \cite{IMCF}. These
minimizing properties characterize how the weak flow differs from the
classical flow. The classical flow runs into trouble when the mean
curvature changes sign. If the mean curvature is negative on a part of
$N_t$, we could decrease the area of $N_t$ be flowing that patch
out. So the weak flow avoids these areas by making sure that $N_t^+$ is
always the outermost surface with its area. This sometimes necessitates
jumping. 

The existence and uniqueness of solutions to equation
\Eqref{IMCF-weak} with initial data are beyond our scope. However, for
completeness, here is the existence and uniqueness theorem (3.1) from
\cite{IMCF}.
\begin{thm}\Label{IMCF-exist}
  Let $M$ be a complete, connected Riemannian $n$-manifold without
  boundary. Suppose there exists a proper, locally Lipschitz, weak
  subsolution of \Eqref{IMCF-IVP} with a precompact initial condition.

Then for any nonempty, precompact, smooth open set $E_0$ in $M$, there
exists a proper, locally Lipschitz solution $u$ of \Eqref{IMCF-IVP}
with initial condition $E_0$, which is unique in $M\setminus
E_0$. Furthermore, the gradient of $u$ satisfies the estimate
\begin{gather}
  \abs{\grad u(x)} \leq \sup_{\p E_0\cap B_r(x)} H_+ + \frac{C(n)}r,
  \quad \text{a.e. }x\in M\setminus E_0, 
\end{gather}
for each $0<r\leq \gs(x)$. 
\end{thm}
The function $\gs(x)$ depends on the Ricci curvature of $M$ near
$E_0$, but is always positive. More importantly, the requirement for a
subsolution is satisfied by any asymptotically flat manifold. A
function like $\ln(R)$ in the asymptotic end will suffice.

\section{Useful Properties of Weak IMCF}
Now that we have outlined the flow, we will describe some useful
properties. 
\begin{lemma}\Label{IMCF-exp-area}
  Let $(E_t)_{t>0}$ solve \Eqref{IMCF-IVP} with initial condition
  $E_0$. As long as $E_t$ remains precompact, we have the following:
  \begin{itemize}
  \item $e^{-t}\abs{\p E_t}$ is constant for $t>0$.
  \item If $E_0$ is its own minimizing hull then $\abs{\p E_t} =
    e^t\abs{\p E_0}$.
  \end{itemize}
\begin{proof}
Since each $E_t$ minimizes the same functional, they must all have the
same value for $J_u(E_t)$. Applying the co-area formula to the
integral in $J_u$ gives
\begin{align}
  J_u(E_t) &= \abs{\partial E_t}- \int_0^t\frac{1}{\abs{\grad u}}
\int_{\p E_s} \abs{\grad u}\de A\de s\\
&= \abs{\partial E_t}- \int_0^t\abs{\p E_s }\de s.
\end{align}
Which has solutions of the form $Ce^t$. By the minimizing hull
properties, we could replace $\p E_t$ with $\p E^+_t$ for $t>0$. Since
$\p E_0^+$ is the limit of $\p E_s$ as $s\searrow 0$, $C=\abs{\p
  E_0^+}$.  If $E_0$ is its own minimizing hull then $\abs{\p
  E_0}=\abs{\p E_0^+}$, and $C=\abs{\p E_0}$.
\end{proof}
\end{lemma}

The next two lemmas tell us that when the classical solution exists,
it agrees with the weak solution.
\begin{lemma}
  Let $(N_t)_{c\leq t<d}$ be a smooth family of surfaces of positive
  mean curvature that solve \Eqref{IMCF-classic} classically. Let
  $u=t$ on $N_t$, $u<c$ inside $N_c$, and $E_t=\{u<t\}$. Then for
  $c\leq t<d$, $E_t$ minimizes $J_u$ in $E_d\setminus E_c$.
\end{lemma}

\begin{lemma}
  Let $E_0$ be a precompact open set in $M$ such that $\p E_0$ is
  smooth with $H>0$ and $E_0=E_0'$. Then any unique solution
  $(E_t)_{0<t<\infty}$ of \Eqref{IMCF-IVP} with initial condition
  $E_0$ coincides with the unique smooth classical solution for a
  short time, provided $E_t$ remains precompact for a short time.
\end{lemma}
The authors of \cite{IMCF} point out that the stopping point for these
theorems will be when either $E_t$ is no longer a minimizing hull, the
mean curvature goes to zero, or the second fundamental form is
unbounded.

The proof of existence and uniqueness of the weak flow are beyond the
scope of this thesis. However, the method is as follows. First they
assume the existence of a subsolution $v$. Then they solve the regularized
problem
\begin{gather}
  \begin{aligned}
    E^\ge u^\ge &:= \Div\left(\frac{\grad u^\ge}{\sqrt{\abs{\grad
            u^\ge}^2 + \ge^2}}\right)-\sqrt{\abs{\grad u^\ge}^2 +
      \ge^2}=0 &\text{in }\gO_L\\
u^\ge &= 0&\text{on }\p E_0\\
u^\ge &= L-2 &\text{on }\p F_L. 
\end{aligned}
\end{gather}
Where $F_L=\{v<L\}$, $\ge$ is small, and $L$ is large, but bounded in
size by a function of $\ge$.  Then they take the limit as
$L\to\infty$, $\ge\to 0$. Some estimates on $\abs{\grad u}$ and $H$
guarantee that the solution passes to the limit.

The authors also proved that in the asymptotic regime, the Hawking
mass of the level sets of the flow converges to the ADM mass of the
manifold

\begin{thm}\Label{ADM-IMCF}
  Assume that $M$ is asymptotically flat and let $(E_t)_{t\geq t_0}$
  be a family of precompact sets weakly solving \Eqref{IMCF-classic}
  in $M$.  Then
\begin{gather}
\lim_{t\to\infty}  m_{\Hk}(N_t) \leq m_\ADM(M).
\end{gather}
\begin{proof}
They show that $N_t$ must approach coordinate spheres in the asymptotic regime. 
This is their Lemma 7.4.
\end{proof}
\end{thm}

\section{Geroch  Monotonicity}
The Geroch Monotonicity formula says that the Hawking mass is
nondecreasing under the inverse mean curvature flow. The original use
of this was to propagate the mass of a surface out to infinity to
compare to the ADM mass using IMCF\@. However, it also provides bounds
on integrals of $H$ and area near the surface.  We will first show how
it arises in the smooth case, and then extend it over the jumps in the
weak flow.

In the smooth case, we simply recall that the Hawking mass is given by 
\begin{gather}
  m_{\Hk}=\sqrt{\frac{\abs{N}}{16\pi}}\left(1-\frac1{16\pi}\int_N
    H^2\right).
\end{gather}
The first variation of area is given by \begin{gather}
\dd{}t \de A = H \gh\de A. 
\end{gather}
Thus under smooth IMCF we have $\dd{}t \de A = \de A$. 
The variation of $H$ is given by
\begin{gather}
  \dd{H}{t} = \gD(-\gh)-\abs{A}^2\gh - \Rc(\gv,\gv)\gh.
\end{gather}
Thus if we look at the integral $\int H^2\de A$ under IMCF we get
\begin{gather}
  \dd{}t \int H^2\de A = \int H^2 
-2\frac{\abs{\grad_N H}^2}{H^2}+2\abs{A}^2
  -2\Rc(\gv,\gv) \de A.
\end{gather}
The Gauss equation contracts to give
\begin{gather}
  K=K_{12}+\gl_1\gl_2=\frac{R}2-\Rc(\gv,\gv)+\frac12(H^2-\abs{A^2})
\end{gather}
Here $K$ and $R$ are the scalar curvatures of $N$ and $M$ respectively,
$K_{12}$ is the sectional curvature of $M$ in the plane tangent to
$N$, and $\gl_i$ are the principal curvatures of $N$.  Using this
equation to cancel the $\Rc$ term gives us the following equation
\begin{align}
  \dd{}t \int_N H^2 &= \int_N 2K-2\frac{\abs{\grad_N H}^2}{H^2}
-\abs{A}^2-R\\
  &= \int_N 4K-R - 2\frac{\abs{\grad_N H}^2}{H^2}
-\frac12\left(\gl_1+\gl_2\right)^2-\frac12\left(\gl_1-\gl_2\right)^2.\\
\intertext{If $R>0$,}
\dd{}t \int_N H^2  &=4\pi \gc(N_t)-\frac12\int_N H^2-\int_N2\frac{\abs{\grad_N
      H}^2}{H^2}+\frac12\left(\gl_1-\gl_2\right)^2.\\
\intertext{If $N_t$ is connected,}
\dd{}t \int_N H^2&\leq 8\pi\left(1-\frac1{16\pi}\int_N H^2\right).
\end{align}
In addition recall that $\abs{N_t}=\abs{N_0}e^t$ in the smooth case.
Thus we get that
\begin{align}
  \sqrt{16\pi}\dd{}t m_{\Hk} &=
  \dd{}t\left[e^{t/2}\left(1-\frac1{16\pi}\int_N
      H^2\right)\right]\\
&\geq\left[\frac12e^{t/2}\left(1-\frac1{16\pi}\int_N
      H^2\right)
-e^{t/2}\frac{8\pi}{16\pi}\left(1-\frac1{16\pi}\int_N H^2\right)
\right]=0.
\end{align}
Hence the Hawking mass is nondecreasing. 

To cover the gap, we simply note that at the jumps, the new surface
$E'_t=E_t^+$ is a minimizing hull for $E_t$. That means that where
their boundaries differ, $\p E'_t$ must have zero mean curvature (else
a variation could decrease its area keeping it outside of $E_t$,
contradicting its minimizing property.) Hence
\begin{gather}
  \int_{\p E_t}H^2\geq \int_{\p E_t'} H^2.
\end{gather}
Thus since $\abs{\p E'_t}=\abs{\p E_t}$ for $t>0$, we see that the,
even at jumps, the Hawking mass can't decrease. This doesn't cover the
possibility of dense jumping, or similar pathological behavior, but
extending the Geroch formula to those cases requires using the
elliptic regularization. 

The statement of Geroch Monotonicity given in \cite{IMCF} for one
boundary component is as follows
\begin{thm}\Label{Geroch}
  Let $M$ be a complete 3-manifold, $E_0$ a precompact open set with
  $C^1$ boundary satisfying $\int_{\p E_0}\abs{A}^2<\infty$, and
  $(E_t)_{t>0}$ a solution to \Eqref{IMCF-IVP} with initial condition
  $E_0$. If $E_0$ is a minimizing hull then
 \begin{multline}
    m_{\Hk}(N_s) \geq m_{\Hk}(N_r) +\\+\frac{1}{(16\pi)^{3/2}}\int_r^s
    \left[16\pi -8\pi \gc(N_t)+ \int_{N_t}\left(2\abs{D\log H}^2 +(\gl_1-\gl_2)^2+R\right)\right]dt
\end{multline}  
for $0\leq r\leq s$ provided $E_s$ is precompact. 
\end{thm}

\chapter{Negative Point Mass Singularity Results}
\Label{IMCF-NPMS}
\section{Negative Point Mass Singularities and IMCF}
Although we do not need this fact, it is interesting to note that near
a negative point mass singularity, we can define the inverse mean
curvature flow. Even though there is no initial surface to start from
we can take a limit of solutions to IMCF for starting surfaces that
converge to $p$.  This actually defines a unique solution $u$ to the
weak inverse mean curvature flow. This was shown in recent work by
Jeffrey Streets \cite{streets}.

\begin{thm}[\cite {streets}] \Label{jeff-exist}
Let $M^3$ be an asymptotically flat manifold with finitely many
singularities at $\{p_i\}$. Then there is a unique solution to
\Eqref{IMCF-weak} on $M\setminus\{p_i\}$.
\end{thm}
In this case since the area of the level sets is exponential in time,
and surfaces near the singularities have vanishing area. Thus this
flow only reaches the singularity at time $-\infty$. 

Streets also showed that these surfaces were the best possible
surfaces, in terms of the Hawking mass, as in the following theorem.
\begin{thm}[\cite{streets}]\Label{jeff-best}
  Let $S_t$ be the family of hypersurfaces defining the solution to
  IMCF above.  Let $P_t$ be any other family of hypersurfaces
  approaching the singularity.  Then, 
  \begin{gather}
	\lim_{t\to-\infty} m_\Hk (P_t) \leq \lim_{t\to-\infty} m_\Hk (S_t) .
  \end{gather}
\end{thm}

We can also extend Geroch Monotonicity down to $t=0$ in the case where
our initial surface has negative Hawking mass.
\begin{lemma}\Label{neg-mass-hull-lemma}
  Let $\gS$ be a surface in an asymptotically flat 3 manifold.  Let
  $\gS'$ be the boundary of the minimizing hull of $\gS$. Let $\gS$ or $\gS'$
  have negative Hawking mass. Then
  \begin{gather}
    m_\Hk(\gS)\leq m_\Hk(\gS').
  \end{gather}

\begin{proof}
  If $\gS'$ has nonnegative Hawking mass then $m_\Hk(\gS')\geq 0 \geq
  m_\Hk(\gS)$ and we are done. Thus we can assume that $\gS'$ has
  negative Hawking mass.  Since $\gS'$ has negative Hawking mass, it
  must intersect $\gS$ on a set of positive measure. Otherwise, $\gS'$
  would be a minimal surface, with Hawking mass
  $\sqrt{\frac{\abs{\gS'}}{16\pi}}>0$.  We define the following sets:
  \begin{gather}
    \gS_0 = \gS'\cap \gS\qquad
    \gS_+ = \gS'\setminus\gS_0\qquad
    \gS_- = \gS\setminus\gS_0
  \end{gather}
  Recalling that $\abs{\gS_+}\leq \abs{\gS_-}$ by the minimization
  property, and that $H=0$ on $\gS_+$, we observe the following:
  \begin{align}
    0>m_\Hk(\gS') &=
 \frac{\sqrt{\abs{\gS_0}+\abs{\gS_+}}}{(16\pi)^{3/2}}
    \left(16\pi-\int_{\gS_0}H^2\right)\\
    &\geq  \frac{\sqrt{\abs{\gS_0}+\abs{\gS_-}}}{(16\pi)^{3/2}}
    \left(16\pi-\int_{\gS_0}H^2\right)\\
    &\geq  \frac{\sqrt{\abs{\gS_0}+\abs{\gS_-}}}{(16\pi)^{3/2}}
    \left(16\pi-\int_{\gS_0}H^2-\int_{\gS_-}H^2\right)=m_\Hk(\gS).
  \end{align}
\end{proof}
\end{lemma}

With this lemma and Geroch Monotonicity we can prove  the following
lemma.
\begin{lemma}\Label{geroch-2}
  Let $(M,g)$ be an asymptotically flat manifold with ADM mass $m$,
  nonnegative scalar curvature and a single regular negative point
  mass singularity $p$. Let $\{\gS_i\}$ be a smooth family of surfaces
  converging to $p$, which eventually have negative Hawking mass. Then
  for sufficiently large $i$,
  \begin{gather}
    m_\Hk(\gS_i)\leq m.
  \end{gather}
  \begin{proof}
    Since for large enough $i$, $\gS_i$ has non-positive Hawking mass, we
    can apply Lemma~\Ref{neg-mass-hull-lemma} to show that $\gS_i'$
    must have larger Hawking mass. From this surface, we start Inverse
    Mean Curvature Flow. Theorem~\Ref{Geroch} tells us that the
    Hawking mass of the surfaces $N_t$ defined by IMCF starting with
    $\gS_i'$ only increase. Theorem~\Ref{ADM-IMCF} tells us that the
    increasing limit of the Hawking masses these surfaces is less than
    the ADM mass. Thus the Hawking mass of the starting surface was
    also less than the ADM mass.
    \end{proof}
\end{lemma}

Now we relate the limit of the Hawking masses to the regular mass.
\begin{lemma}\Label{hawking-regular}
  Let $(M,g)$ be an asymptotically flat manifold with nonnegative
  scalar curvature and a single regular negative point mass
  singularity $p$.  Then there is a smooth family of surfaces
  $\{\gS_i\}$ converging to $p$ such that
\begin{gather}
  \lim_{i\to\infty} m_\Hk(\gS_i)=
  -\frac14\left(\frac1\pi\int_{\bar\gS} \bar\gv(\bar\gf)^{4/3}\bar{\de
      A}\right)^{3/2}=m_{\NPMS}(p).
\end{gather}
\begin{proof}
  The
Hawking mass of a surface $\gS_i$ is given by
\begin{gather}
  m_\Hk(\gS_i) =
  \sqrt{\frac{\abs{\gS_i}}{16\pi}}\left(1-\frac1{16\pi}\int_{\gS_i}H^2\de
    A\right).
\end{gather}
Since the areas of the surfaces are converging to zero  we have
\begin{gather}
  \lim_{i\to\infty}  m_\Hk(\gS_i) =
-  \lim_{i\to\infty} \frac{\sqrt{\abs{\gS_i}}}{(16\pi)^{3/2}}\int_{\gS_i}H^2\de A.
\end{gather}
By the H\"older inequality this is bounded as follows
\begin{gather}\Label{loc-CS}
  - \frac{\sqrt{\abs{\gS_i}}}{(16\pi)^{3/2}}\int_{\gS_i}H^2\de A\leq -
  \frac1{(16\pi)^{3/2}}\left(\int_{\gS_i} H^{4/3}
    \de A\right)^{3/2}.
\end{gather}
Switching to the resolution space, we use the formula
\begin{gather}
  H =\bar\gf^{-2} \bar H + 4\bar\gf^{-3}\bar\gv(\bar\gf ).
\end{gather}
Putting this into the previous equation we get
\begin{align}
   \int_{\gS_i} H^{4/3}
    \de A&=\int_{\bar{\gS}_i} \left( \bar\gf^{-2} \bar H +4
      \bar\gf^{-3}\bar\gv(\bar\gf ) \right)^{4/3}\gf^4\bar{\de A}\\
    &=\int_{\bar{\gS}_i} \left( \bar\gf \bar H +4
      \bar\gv(\bar\gf ) \right)^{4/3}\bar{\de A}.
\end{align}
Since $\bar\gf$ is zero on $\bar\gS$ and $\bar H$ is bounded, the first
term goes to zero. The second term converges since the family of
surfaces $\{\gS_i\}$ are converging smoothly.
\begin{gather}
  \lim_{i\to\infty}\int_{\bar{\gS}_i} \left( \bar\gf \bar H +4
    \bar\gv(\bar\gf ) \right)^{4/3}\bar{\de A} =4^{4/3}\int_{\bar\gS}
  \bar\gv(\bar\gf )^{4/3}\bar{\de A}.
\end{gather}
Combining all of these equations we have
\begin{gather}
  \lim_{i\to\infty} m_\Hk(\gS_i)\leq
  -\frac14\left(\frac1\pi\int_{\bar\gS} \bar\gv(\bar\gf)^{4/3}\bar{\de
      A}\right)^{3/2}=m_{\NPMS}(p).\Label{loc-HP}
\end{gather}
To see when this estimate is sharp, we look at inequality
\Eqref{loc-CS} since that is the only inequality is our estimate. In
the limit, this inequality is an equality exactly when the ratio of
the maximum and minimum values of $H$ approaches $1$. We choose a
resolution such that $\bar\gv(\bar\gf)=1$ on the boundary. We also
choose a family of surfaces $\gS_i$ given by level sets of $\bar\gf$.
Then if we look at the ratio
\begin{gather}
  \lim_{\gf\to0} \frac{H_\text{min}}{H_{\text{max}}}=
  \lim_{\gf\to0}\frac{\bar\gf \bar
    H_\text{min}+4\bar\gv(\bar\gf)}{\bar\gf \bar
    H_\text{max}+4\bar\gv(\bar\gf)},
\end{gather}
and remember that $\bar H$ is bounded, we see that the
$\bar\gv(\bar\gf)$ terms dominate, and as $\bar\gf\to0$, this ratio
approaches 1. Thus with this resolution and this family of surfaces,
inequality \Eqref{loc-HP} will turn to an equality.
\end{proof}
\end{lemma}
With these results we can prove
the following theorem
\begin{thm}\Label{penrose-one-regular}
  Let $(M,g)$ be an asymptotically flat manifold with nonnegative
  scalar curvature and a single regular negative point mass
  singularity $p$.  Then the ADM mass of $M$ is at least the mass of
  $p$.
\begin{proof}
  First consider the case when $p$ can be enclosed by a surface,
  $\gS$, with nonnegative Hawking mass. The minimizing hull of a
  surface with nonnegative Hawking mass has nonnegative Hawking mass.
  Thus we can run IMCF from $\gS'$, and the AMD mass of $M$ is at
  least $m_\Hk(\gS')\geq0$.  However, the regular mass of $p$ is
  always nonpositive so in this case we are done.

  Now assume that $p$ cannot be enclosed by a surface with nonnegative
  Hawking mass. By Lemma~\Ref{geroch-2} we know that the ADM mass is
  greater than the Hawking masses of any sequence of surface
  converging to $p$ which have negative Hawking mass. By
  Lemma~\Ref{hawking-regular} we know that there is a family of
  surfaces converging to $p$ which have the mass of $p$ as the limit
  of their Hawking mass, hence the ADM mass is greater then their
  Hawking masses which limit to the regular mass.
\end{proof}
\end{thm}

This can be extended to general negative point mass
singularities. However, first we need to consider the effect of
multiplication by a harmonic conformal factor on the ADM mass of a
manifold. 
\begin{lemma}\Label{harm-mass-mod}
  Let $(M^3,g)$ be an asymptotically flat manifold. Let $\gf$ be a
  harmonic function with respect to $g$ with asymptotic expansion
  \begin{gather}
    \gf = 1+\frac{C}{\abs{x}_g}+\bigo\left(\frac{1}{\abs{x}^2_g}\right).
  \end{gather}
  Then, if the ADM mass of
  $(M^3,g)$ is $m$, the ADM mass of $(M^3,\gf^4g)$ is $m+2C$. 
  \begin{proof}
    This is a direct calculation. We write $g^\gf = g \gf^4$, and
    calculate, only keeping the terms of lowest order in $1/\abs{x}$
    since we are taking limits as $\abs{x}\to\infty$.
\begin{align}
  m_\gf&=
  \lim_{\abs{x}\to\infty}\frac1{16\pi}\int_{S^\gd}\left(g_{ij,i}^\gf
    -g_{ii,j}^\gf\right)n^j\,dA\\
  &=
  \lim_{\abs{x}\to\infty}\frac{\gf^4}{16\pi}\int_{S^\gd}\left(g_{ij,i}
    -g_{ii,j}\right)n^j\,dA+
  \lim_{\abs{x}\to\infty}\frac{\gf^3}{4\pi}\int_{S^\gd}\left(\gd_{ij}\gf_i
    -\gd_{ii}\gf_j\right)n^j\,dA\\
  &= \lim_{\abs{x}\to\infty}\gf^4m +
  \lim_{\abs{x}\to\infty}\gf^3\lim_{\abs{x}\to\infty}
  \frac{1}{4\pi}\int_{S^\gd}\left(\gf_j-3\gf_j\right)n^j\,dA\\
  &=m +
  \lim_{\abs{x}\to\infty}\frac{1}{2\pi}\int_{S^\gd}\gf_jn^j\,dA\\
  &=m + \lim_{\abs{x}\to\infty}\frac{1}{2\pi}\int_{S^\gd}\ip{\grad\gf}{\gv}\,dA\\
  &=m+2C.
\end{align}

\end{proof}
\end{lemma}

Using this we can now extend Theorem \Ref{penrose-one-regular} to a
general singularity. 
\begin{thm}\Label{penrose-one-general}
  Let $(M,g)$ be an asymptotically flat manifold with nonnegative
  scalar curvature and a single negative point mass singularity $p$.
  Then $m$, the ADM mass of $M$, is at least the mass of $p$.
  \begin{proof}
    If the capacity of $p$ is nonzero, then the statement is trivial.
    Thus we assume the capacity of $p$ is zero.  Using the terminology
    of Definition \Ref{NPMS-mass-surf}, Theorem
    \Ref{penrose-one-regular} tells us that the ADM mass of
    $(M,h^4_ig)$ is at least the mass of the regular singularity at
    $\gS_i=p_i$. Each $h_i$ is defined by the equations
    \begin{align}
      \gD h_i &= 0 \\
      \lim_{x\to\infty}h_i &=1\\
      h_i &= 0 \text{ on } \gS_i .
    \end{align}
    Thus it has asymptotic expansion
    \begin{gather}
      h_i = 1-\frac{C_i}{\abs{x}}+\bigo\left(\frac1{\abs{x}^2}\right).
    \end{gather}
    Where $C_i$ is the capacity of $\gS_i$.  Thus, the ADM mass,
    $m_i$, of $(M,h^4_ig)$ is given by $m-2C_i$. Now we know that
    $m_i\geq m_{\NPMS}(p_i)$. Taking $\varlimsup$ of both sides gives
    us
    \begin{gather}
      \varlimsup_{i\to\infty}m_i\geq \varlimsup_{i\to\infty} m_{\NPMS}(p_i)
    \end{gather}
    Since $C_i$ is going to zero, the left hand side is simply $m$,
    and so has no dependence on which $\{\gS_i\}$ we chose in our mass
    calculation. Thus we get
    \begin{gather}
      m \geq \sup_{\{\gS_i\}}\varlimsup_{i\to \infty} m_{\NPMS}(p_i).
    \end{gather}
    as desired.
  \end{proof}
\end{thm}

\section{Capacity Theorem}

Perhaps the most important new results of this thesis are Theorems
\Ref{cap-mass-lemma} and \Ref{cap-mass-theorem} which relate the
capacity of a point to the Hawking masses of surfaces near that
point.

The capacity of a surface provides a measure of its size as seen from
infinity. We extend the definition of the capacity of surface to the
capacity of a negative point mass singularity.  We then show that if a
NPMS has non-zero capacity the Hawking mass of any family of surfaces
converging to it must go to negative infinity.  First the definition
of capacity:
\begin{defn}
  Let $\gS$ be surface in an asymptotically flat manifold
  $M$. \fluff{Choose an end of $M$.}Define the \emph{capacity} of $\gS$ by
  \begin{gather}
C(\gS)= \inf\left\{\left. \int_M \norm{\nabla \gf}^2 dV \right| 
        \gf(\gS)=1, \gf(\infty)=0\right\}.
\end{gather}
\end{defn}
It is worth noting that if $\gS$ and $\gS'$ are two surfaces in $M$ so that
$\gS$ divides $M$ into two components, one containing infinity
and the other containing $\gS'$, then
\begin{gather}
  C(\gS')\leq C(\gS)
\end{gather}
since the infimum is over a larger set of functions.

We will next define the capacity of a singular point. The natural
definition is the one we want.
\begin{defn}
  Let $p$ be singular point in an asymptotically flat manifold $M$.
  \fluff{Chose an end of $M$.} Chose a sequence of surfaces $\gS_i$ of
  decreasing diameter enclosing $p$. Then define the \emph{capacity}
  of $p$ by the limit of the capacities of $\gS_i$. \fluff{both with respect
  to the chosen end.}
\end{defn}
Before using this definition we have to show that it is unique.
\begin{lemma}
  Let $\gS_i$ and $\widetilde\gS_i$ be two sequences of surfaces
  approaching the point $p$. If $\lim C(\gS_i)=K$, $\lim
  C(\widetilde\gS_i)=K$. Hence $C(p)$ is well defined.
\begin{proof}
  Since the $\gS_i$ are going to $p$, for any given $\widetilde
  \gS_{i_0}$, we can choose $i_0$ such that for all $i> i_0$, $\gS_i$
  is contained within $\widetilde \gS_{i_0}$. Thus if $\gf$ is a
  capacity test function for $\widetilde \gS_{i_0}$, i.e.\
  $\gf(\gS_{i_0})=1$ and $\gf\to0$ at infinity, then $\gf$ is also a
  capacity test function for $\gS_i$. Since $C(\gS_i)$ is taking the
  infimum over a larger set of test functions than $C(\widetilde
  \gS_{i_0})$, $C(\gS_i)\leq C(\widetilde \gS_{i_0})$. Thus if we
  create a new sequence of surfaces $\bar \gS_i$, alternately choosing
  from $\gS_i$ and $\widetilde\gS_i$, such that each surface contains
  the next we get a nonincreasing sequence of capacities. Thus if
  either original sequence of surfaces has a limit of capacity, then
  this new sequence must as well, and it must be the same. Hence,
  $\lim_{i\to\infty}C(\gS_i)=\lim_{i\to\infty}C(\widetilde\gS_i)$.
\end{proof}
\end{lemma}

Now we look at the relationship between capacity and the Hawking mass
of a surface. We will use techniques similar to those used in
\cite{Bray-Miao}.

\begin{thm}
  \Label{cap-mass-lemma}
  Let $M$ be an asymptotically flat 3 manifold with nonnegative scalar
  curvature, and negative point mass singularity $p$.  Let $\gS_i$ be
  a family of surfaces converging to $p$.  Assume each $\gS_i$ is a
  minimizing hull. Assume the areas of $\gS_i$ are going to zero. Then
  if the Hawking mass of the surfaces is bounded below, the capacities
  of surfaces of foliation near $p$ must go to zero.
  \begin{proof}
    To use Geroch monotonicity, we need to know that our IMCF surfaces
    stay connected. In the weak formulation of IMCF, the level sets
    $\gS_t$ always bound a region in $\bar M$.  Thus if $\gS_t$ is not
    connected, one of its components $\gS_t^*$ must not bound a
    region. That is, $\gS_t^*$ is not homotopic to a point in $M$.
    Since $M$ is smooth, it must have finite topology on any bounded
    set. Thus we know that near $p$, there is a minimum size for a
    surface that does not bound a region. Call this size
    $A_\text{min}$.  Thus if we have any surface that does not bound a
    region, it must have area greater then $A_{\text{min}}$. According
    to Lemma \Ref{IMCF-exp-area} the area of our surfaces grow
    exponentially. Thus if we restrict ourselves to starting IMCF with
    a surface with area $A_\text{min}/e$, and only run the flow for
    time $1$, we will stay connected. At first glance it seems we may
    need to worry about the jumps in weak IMCF, however Geroch
    monotonicity doesn't depend on smoothness of the flow, and neither
    does the area growth formula.  Lemma \Ref{IMCF-exp-area} holds
    from $t=0$.  Thus even with jumps, the area of our surfaces will
    remain below $A_\text{min}$.
    
    Now recall that capacity of a surface is defined by
    \begin{gather}
      C(\gS)= \inf\left\{\left. \int_M \norm{\nabla \gf}^2 dV \right| 
        \gf(\gS)=1, \gf(\infty)=0\right\}.
    \end{gather}
    Here, the integral is only over the portion of $M$ outside of
    $\gS$.  Call this integral, $\E(\gf)$, the \emph{energy} of $\gf$.
    Thus for any $\gf$ with $\gf(\infty)=0$ and $\gf(\gS)=1$ we have
    $\E(\gf)\geq C(\gS)$. So we will find an estimate that relates the
    Hawking mass and the energy of a test function $\gf$.
    
    Choose a starting surface $\gS$ with sufficiently small starting
    area. Let $f$ be the level set function of the associated weak
    IMCF starting with the surface $\gS$. Call the resulting level
    sets $\gS_t$. Now if we use a test function of the form
    $\gf=u(f)$, then the energy of $\gf$ is given by
    \begin{gather}
      \E(\gf)= \int_M \norm{\N f}^2(u')^2 dV.
    \end{gather}
    Since $f$ is given by IMCF, we know that $\norm{\N f} = H$ where
    $H$ is the mean curvature of the level sets. Next we will use the
    co-area formula. This states that if we have a function $z$ on a
    domain $\gO$, and a function $h:\Re\to\Re$ so that the range of
    $h(z)$ is $[a,b]$, then
    \begin{gather}
      \int_\gO h \,dV = \int_a^b h(z(t))\int_{S_{t}}\abs{\grad z(t)}\,dA_t\,dt.
    \end{gather}
    Here $S_t$ are the level sets of $h(z(t))$.  If we use the
    co-area formula with the foliation $\gS_t$, our integral becomes
    \begin{gather}
      \E(\gf)= \int_0^\infty (u'(t))^2\int_{\gS_t} \abs{H} dA_t\,dt.
    \end{gather}
    Here the co-area gradient term cancels one of the $\norm{\N f}$ terms.
    Now we will bound the interior integral of curvature. We know that
    IMCF causes the Hawking mass to be nondecreasing in $t$. We first
    rewrite the definition of the Hawking mass $m_\Hk(\gS_t^i) = m(t)$
    as:
    \begin{gather}
      \int H^2 dA_t = 16\pi\left(1 - m(t)\sqrt{\frac{16\pi}{A(t)}}\right).
    \end{gather}
    Here $A(t)$ is the area of $\gS_t$.  Since the Hawking mass is
    nondecreasing under IMCF we have:
    \begin{gather}
      \int H^2 dA_t \leq 16\pi\left(1 - m(0)\sqrt{\frac{16\pi}{A(t)}}\right).
    \end{gather}
    Thus we can use Cauchy-Schwartz to get:
    \begin{gather}
      \int \abs H dA_t \leq \sqrt{A(t)}\sqrt{ 16\pi\left(1 -
          m(0)\sqrt{\frac{16\pi}{A(t) }}\right)}.
    \end{gather}
    We can rewrite this as:
    \begin{gather}
      \int \abs H dA_t \leq \sqrt{\ga A(t)+\gb\sqrt{A(t)}}.
    \end{gather}
    Furthermore, since $A(t)$ grows exponentially in $t$, we can write
    this as:
    \begin{gather}
      \int \abs H dA_t \leq \sqrt{\ga e^t+\gb e^{t/2}}=v(t).
    \end{gather}
    Where $A_0$ has been absorbed into $\ga$ and $\gb$. Thus our
    energy formula has become
    \begin{gather}
      \E(\gf) \leq \int_0^\infty (u'(t))^2v(t)\,dt.
    \end{gather}
    with
    \begin{gather}
      v(t)=\sqrt{\ga e^t+\gb e^{t/2}}
    \end{gather}
    where $\ga= 16\pi A_0$, $\gb = (16\pi)^{3/2}A_0^{1/2}\abs{m_0}$,
    and $A_0$ is $A(\gS_0)$. This means we can pick our test
    function $u(t)$ to be as simple as:
    \begin{gather}
      u(t)=\begin{cases} 1-t & 0\leq t\leq 1\\ 0 &t\geq 1
      \end{cases}
    \end{gather}
    Then our integral becomes:
    \begin{align*}
      E(\gf)&\leq \int_0^1 v(t)dt\\
      &=\int_0^1 \sqrt{\ga e^t+\gb e^{t/2}}\,dt\\
      &=\int_0^1 e^{t/4}\sqrt{\ga e^{t/2}+\gb}\,dt\\
      &=4\int_{1}^{e^{1/4}} \sqrt{\ga x^2+\gb}\,dx\\
\fluff{      &=2\left.\left.x\sqrt{\ga x^2+\gb}+
          \frac{\gb\asinh\left(\sqrt{\frac\ga\gb}x\right)}
          {\sqrt\ga}\right.\right|_{1}^{e^{1/4}}.}
    &\leq 4\int_1^{e^{1/4}}\sqrt\ga x +\sqrt\gb\\
    &=2\sqrt{\ga}(e^{1/2}-1)+4\sqrt\gb(e^{1/4}-1)\\
    &\leq 2\sqrt\ga+2\sqrt\gb.
    \end{align*}
    Since $m_{\Hk}(\gS)\leq \sqrt{\frac{\abs{\gS}}{16\pi}}$, $m_0$ is
    bounded above. By assumption $m_0$ is bounded below, so $\ga$ and
    $\gb$ are bounded by multiples of $A_0$ and $\sqrt{A_0}$
    respectively. \fluff{  Thus as $A_0\to0$, the only term that could
    possible cause trouble, the last term, goes like 
    \begin{align*}
      \frac{\gb\asinh\left(\sqrt{\frac\ga\gb}x\right)}
      {\sqrt\ga}
      &=\frac{(16\pi)^{3/2}A_0^{1/2}\abs{m_0}
        \asinh\left(x\sqrt{\frac{16\pi A_0}{(16\pi)^{3/2}A_0^{1/2}
              \abs{m_0}}}\right)}{\sqrt{16\pi A_0}}\\
      &=16\pi\abs{m_0}
      \asinh\left(x\frac{A^{1/4}_0}{(16\pi)^{1/4}\abs{m_0}^{1/2}}\right)\\
      & \leq  16\pi\abs{m_0}
      \left(x\frac{A^{1/4}_0}{(16\pi)^{1/4}\abs{m_0}^{1/2}}\right)\\
      &=(16\pi)^{3/4}\abs{m_0}^{1/2}A_0^{1/4}x.
    \end{align*}
    Combining this with the other term we get
    \begin{align*}
      C(\gS)&\leq \E(\gf) \\
      &\leq 2e^{1/4}\sqrt{16\pi A_0 e^{1/2}+ (16\pi)^{3/2}A_0^{1/2}\abs{m_0}}+
      (16\pi)^{3/4}\abs{m_0}^{1/2}A_0^{1/4}e^{1/4}.
    \end{align*}} 
  Thus $\E(\gf)$ goes to zero if $A_0\to 0$ and $m_0$
  is bounded.  Hence $C(p)$ must be zero since it is the infimum over
  a positive set with elements approaching zero.
  \end{proof}
  
\end{thm}

\begin{thm}[Capacity Theorem]
  \Label{cap-mass-theorem}
  Let $M$ be an asymptotically flat 3 manifold with nonnegative scalar
  curvature, and negative point mass singularity $p$, such that there
  exists a family of surfaces, $\gS_i$, converging to $p$ with area
  going to zero.  Then if the capacity of $p$ is nonzero, the Hawking
  masses of the surfaces $\gS_i$ must go to $-\infty$.
  \begin{proof}
    Any such family of surfaces will generate a family $\gS_i'$ of
    minimizing hulls that will also converge to $p$.  By Lemma
    \Ref{cap-mass-lemma}, the masses of $\{\gS_i'\}$ must go to
    $-\infty$. Thus the masses of $\{\gS_i'\}$ must go to $-\infty$.
    Thus for sufficiently large $i$, the masses of the minimizing
    hulls are all negative.  From then on Lemma
    \Ref{neg-mass-hull-lemma} applies, and the masses of $\gS_i$ must
    be less then the masses of $\gS_i'$. Hence they also converge to
    $-\infty$.
  \end{proof}

\end{thm}

\chapter{Symmetric Singularities}
\Label{Axi}
In this chapter we will look at some more examples of negative point
mass singularities. The symmetry ansatz provides more structure
than in the definition of a NPMS.  

\section{Spherical Solutions}
A spherically symmetric manifold, $M$, has a metric given by
\begin{gather}
  ds^2 = dr^2+\frac{A(r)}{4\pi}dS^2.
\end{gather}
We can directly calculate that the scalar curvature of $M$ is given by
\begin{gather}
R= \frac{16\pi A+A'^2-4AA''}{2A^2}.
\end{gather}
For this manifold to be asymptotically flat it is necessary for the
Hawking masses of the coordinate spheres to approach a constant. The
Hawking mass of a coordinate sphere is given by
\begin{gather}
  m_\Hk(S)=\sqrt{\frac{A}{16\pi}}\left(1-\frac{1}{16\pi}\frac{A'^2}{A}\right).
\end{gather}
Due to the spherical symmetry of the manifold, we know that if IMCF is
started with coordinate spheres, it must continue with coordinate
spheres. Thus we know that this quantity must be non-decreasing.  For
this manifold to be asymptotically flat, this quantity must have a
limit at $\infty$. This limit is the ADM mass.

The only possible location for a singularity in such a manifold is at
the origin. We can find a straight forward function to resolve the
singularity. We need a smooth function $\gf$ such that
\begin{gather}
\lim_{r\to\infty}  \frac{A(r)}{\gf^4}=4\pi,\Label{sym-res-def}
\end{gather}
or any other constant. 

Multiplying by $\gf^{-4}$ to find the model space changes our metric
to the form
\begin{gather}
  \tilde ds^2 = d\gr^2 + \frac{\tilde A(\gr)}{4\pi} dS^2.
\end{gather}
Where $\tilde A$ goes to $4\pi$ as $\gr$ approaches zero. Let $\gS$ be
the surface $\gr=0$.  The behavior of $\gf$ away from the singularity
is very flexible as long it is bounded, nonzero, and goes to one at
infinity. In order for this to be a regular singularity, we need that
$\gr$ be well defined. Thus we must require
\begin{gather}
  \gr(r) = \int_0^r\gf^{-2}dr = \int_0^r\frac{dr}{A^{1/2}(r)}<\infty
\end{gather}
for finite $r$. Any spherically symmetric singularity with this
condition on $A(r)$ must be regular. The regular mass is given by
\begin{gather}
  m_\NPMS = -\frac14\left(\frac1\pi\int_{\gS} \tilde
    \gv(\gf)^{4/3}d\tilde A \right)^{3/2}.
\end{gather}
Since we have only defined the asymptotic behavior of $\gf$ near
$r=0$, we will compute everything as limits as $r$ goes to zero. First
we need to find $\tilde\gv(\gf)$. Since $\gr= \int_0^r\gf^{-2}dr$, a
chain rule calculation of $\pp{\gf}{\gr}$ yields
\begin{gather}
  \tilde\gv(\gf) = \gf^2 \pp{\gf}{r}.
\end{gather}
Putting this in and expanding the other terms in the definition of the
regular mass gives us
\begin{align}
  m_\NPMS &=-\lim_{r\to0}  \frac14\left(\frac1\pi\int_{\gS_r}
\left(\gf^{2}\pp{\gf}{r}\right)^{4/3}\gf^{-4}dA \right).^{3/2}\\
&=-\lim_{r\to0}  \frac14\frac{1}{\pi^{3/2}}
A^{3/2}(r)\gf^{-2}\left(\pp{\gf}{r}\right)^2
\end{align}
Using equation \Eqref{sym-res-def} and l'H\^opital's rule we find that
\begin{gather}
  \pp{\gf}{r} =\frac{1}{16\pi\gf^3}\pp{A}{r} 
\end{gather}
Continuing from above
\begin{gather}\Label{sphere-mass}\begin{aligned}
   m_\NPMS &= -\lim_{r\to0} \frac14\frac{1}{\pi^{3/2}}
\frac{1}{256\pi^2}\gf^{-8}A^{3/2}A'^2\\
 &=-\lim_{r\to0}\frac{1}{1024\pi^{7/2}}
\frac{16\pi^2}{A^2}A^{3/2}A'^2\\
&=-\lim_{r\to0}\frac1{64\pi^{3/2}}\frac{A'^2}{A^{1/2}}.
\end{aligned}\end{gather}
This agrees with the limit of the Hawking masses of coordinate
spheres as $r\to0$ and hence $A\to 0$ as well.

For completeness, we can also examine the capacity of the central
point in these solutions. The capacity of a coordinate sphere in an
asymptotically flat spherically symmetric manifold is given by the
$1/\gr$ term in of the harmonic function that is $1$ on the sphere and
$0$ at infinity. Since harmonic functions have constant flux with
respect to the outward normal, we could also describe this function as
the constant flux function which goes to $0$ at infinity and $1$ on
the coordinate sphere. Then the capacity of the sphere is given by the
flux constant of this function, divided by $-4\pi$. Reversing this
definition we define the following function:
\begin{defn}
  Let $(M,g)$ be an asymptotically flat spherically symmetric
  manifold. Let $f$ be the radial function that has constant outward
  flux $-4\pi$ through coordinate spheres and goes to zero at
  infinity.  Call $f$ the \emph{radial capacity function} for $(M,g)$.
\end{defn}
This definition allows the following lemma.
\begin{lemma}
  Let $f(r)$ be the radial capacity function for the manifold
  $(M,g)$.  Then the capacity of the coordinate sphere at $r=r_0$ is
  given by $f(r_0)^{-1}$.
  \begin{proof}
    Consider the function $f(r)/f(r_0)$. This function is $1$ on the
    coordinate sphere $r=r_0$, goes to zero at infinity, and is
    harmonic. In the asymptotic regime its $1/\gr$ coefficient is
    $1/f(r_0)$. Thus this is the capacity of the sphere $r=r_0$.
  \end{proof}
\end{lemma}
Thus the capacity of the central point is given by $\lim_{r\to
  0}f(r)^{-1}$. To calculate this value we first note that
\begin{gather}
  \dd{f}{r} = -\frac{4\pi}{A(r)}
\end{gather}
since $f$ has constant flux $-4\pi$. Thus 
\begin{gather}
  f(r) = -\int_r^{\infty} f'(r) dr = 4\pi\int_r^\infty \frac{dr}{A(r)}.
\end{gather}
In order for the capacity of the central point to be nonzero, this
must be finite. Asymptotic flatness tells us that the part of the
integral in the asymptotic regime is finite. Thus the only concern is
where $r\to0$ and hence $A(r)\to 0$. Thus if
\begin{gather}
  \int_0^\ge \frac{dr}{A(r)}
\end{gather}
is finite, our central point has positive capacity. For example, if we
assume that as $r\to0$, $A(r)$ is asymptotically a multiple of a power
of $r$, as $kr^p$. Then the capacity is positive exactly when $p<1$.
In this case, we see by equation~\Eqref{sphere-mass} that the mass of
the singularity is infinite, confirming
Theorem~\Ref{cap-mass-theorem}. However if $1\leq r<4/3$, we see that
the mass of our singularity is still infinite, but the capacity is now
finite. This removes the possibility of strengthening
Theorem~\Ref{cap-mass-theorem} into an if and only if without additional
hypotheses.

\section{Overview of Weyl Solution}
We will be looking at axisymmetric static vacuum spacetimes. These
examples have two Killing fields, one spacelike and one timelike.
These reflect the rotational and time translation symmetry of our
spacetime. Furthermore the timelike Killing field is hypersurface
orthogonal, and the two Killing fields commute. Since our two vector
fields commute, we can call the timelike field $\p_t$ and the
spacelike field $\p_\gth$.  Since they are Killing fields, we know the
metric is only a function of the remaining two coordinates $x_1,x_2$. 
We can also assume that $\p_\gth$ is orthogonal to $\p_t$. Thus our
metric is of the form
\begin{gather}
 g= -A^2dt^2+B^2d\gth^2+g_{11}dx_1^2+g_{12}dx_1dx_2+g_{22}dx_2^2.
\end{gather}
Here $A,B,g_{ij}$ are functions of $x_1$ and $x_2$. We set
$x_1=AB$, and chose $x_2$ to be orthogonal to the other
coordinates. Renaming $x_1$, $\gr$ and $x_2$, $z$, our metric takes
the form
\begin{gather}
  g=-A^2dt^2+\gr^2A^{-2}d\gth^2+U^2d\gr^2+V^2dz^2.
\end{gather}
See Theorem 7.1.1 in \Cite{wald} for an explanation of the lack of
cross terms and further details on this derivation. 

We may also rescale $z$ to set $V=U$. Our metric is now encoded in the
functions $A(\gr,z)$ and $U(\gr,z)$. If we define $\gl$ and $\gm$ by
\begin{gather}
\gl = \ln{A}\qquad \gm = \ln(AU)
\end{gather}
our metric looks like
\begin{gather}
  g=-e^{2\gl}dt^2+e^{-2\gl}\left[\gr^2d\gth^2+
e^{2\gm}\left(d\gr^2+dz^2\right)\right].
\end{gather}
Now, if we compute the curvature of this metric, and apply the vacuum
condition we get the following equations for $\gl$ and $\gm$:
  \begin{align}
0&=  \gl_{\gr\gr}+\frac1\gr\gl_\gr+\gl_{zz}\Label{gl-harm}\\
\gm_\gr &= \gr\left(\gl_\gr^2-\gl_z^2\right)\Label{gm-z}\\
\gm_z & = 2\gr\gl_\gr\gl_z\Label{gm-gr}.
  \end{align}
  The first one is the same as the statement that $\gD\gl=0$ when
  viewed as a function of flat $\Re^3$ with coordinates $(r,z,\gth)$.
  We can use this to generate spacetimes. We can think of $\gl$ as a
  classical potential function. However, we should not think that $g$
  gives the metric with this gravitational potential. For example, the
  Schwarzschild solution is generated by a $\gl$ that is not
  spherically symmetric

  As always we are interested in asymptotically flat spacetimes. In
  this case we want our metric to be the flat metric in cylindrical
  coordinates at infinity. Thus we need the following asymptotics on
  $\gl$ and $\gm$:
\begin{align}
  \lim_{r\to\infty} \gl &= 0\\
  \lim_{r\to\infty} \gm &= 0
\end{align}
Where $r^2=\gr^2+z^2$. Since $\gl$ is flat-harmonic, we know that at
infinity it looks like $C\abs{r}^{-1}+\bigo(r^{-2})$.  This
decay and equation~\Eqref{gl-harm} tells us that
\begin{align}
  \lim_{r\to\infty} \abs{\gl_{\gr}},\abs{\gl_z}\leq \frac{C}{r^2}.
\end{align}
Putting those into equations \Eqref{gm-z} and \Eqref{gm-gr} gives us
\begin{align}
  \lim_{r\to\infty} \abs{\gm_{\gr}},\abs{\gm_z}\leq \frac{C}{r^3}.
\end{align}
These conditions are enough for asymptotic flatness. Thus all that is
required of our metric  metric for it to be asymptotically flat
is that $\gl$ and $\gm$ approach zero at $\infty$. To calculate the ADM
mass of $g$, we need coordinates that are asymptotically Cartesian
rather then asymptotically cylindrical. Using the change of
coordinates
\begin{gather}
  x=\gr\cos\gth \qquad y = \gr\sin\gth
\end{gather}
we get that our metric is
\begin{multline}
g =-e^{2\gl}dt^2
+e^{-2\gl}\left(1+\frac{\left(e^{2\gm}-1\right)x^2}{x^2+y^2}\right)dx^2
+ \frac{e^{-2\gl}\left(e^{2\gm}-1\right) x y}{x^2+y^2}dx\,dy\\
+e^{-2\gl}\left(1+\frac{\left(e^{2\gm}-1\right)y^2}{x^2+y^2}\right)dy^2
+e^{-2\gl}e^{2\gm}dz^2
\end{multline}
If we plug this into the formula for the ADM mass we get the following:
\begin{gather}
  m=\lim_{r\to\infty}\int_{S_r}\frac{e^{-2\gl}}{16\pi r}\left[1-2z\gl_z-2\gr\gl_\gr
+e^{2\gm}\left(2\gr\gm_\gr+2z\gm_z-2\gr\gl_\gr-2z\gl_z
-1\right)\right]\,dA_\gd
\end{gather}
As $r$ grows the terms $\gl_z$ and $\gl_\gr$ are at most order
$r^{-2}$. The derivatives of $\gm$ are at most order $r^{-3}$. The
function $e^{-2\gl}$ is 1, as is the function $e^{2\gm}$. Therefore this
integral is 
\begin{gather}
m
=\frac{-1}{4\pi}\lim_{r\to\infty}\int_{S_r}\frac1r\left(z\gl_z+\gr \gl_\gr\right)\,dA_\gd
=\frac{-1}{4\pi}\lim_{r\to\infty}\int_{S_r}\left<\gv,\nabla \gl\right>_\gd dA_\gd.
\end{gather}
Since $\gl$ is harmonic in the flat metric, we can compute this
on any surface homotopic to a large sphere at infinity.
\begin{lemma}
  Let $(M,g)$ be the $t=0$ slice of an asymptotically flat
  axisymmetric vacuum static manifold with metric
  \begin{gather}
    g=e^{-2\gl}\left[\gr^2d\gth^2+
e^{2\gm}\left(d\gr^2+dz^2\right)\right].
  \end{gather}
Then the $ADM$ mass of $(M,g)$ is given by
\begin{gather}
 -\frac{1}{4\pi}\lim_{r\to\infty}\int_{\gS}\left<\gv,\nabla \gl\right>_\gd dA_\gd,
\end{gather}
where $\gS$ is any surface enclosing the singularities of $\gl$.
\end{lemma}

As in the spherically symmetric case, a useful harmonic function can
tell us about the capacity of the singular points. In this case our
metric comes with a harmonic function. The function $e^\gl$ is
harmonic in our metric.  It goes to $1$ at infinity and infinity at
any positive singularities of our potential. Thus the capacity of
level sets of $e^\gl$ must go to zero as $\gl\to\infty$. This tells us
that these singularities have zero capacity.

\section{Zippoy--Voorhees Metrics}
\fluff{
As mentioned above, the function $\gl$ can be thought of as a
Newtonian potential for $g$. Thus any cylindrically symmetric mass
distribution will give rise to a Weyl metric. However, if we want our
metric to be vacuum, the mass distribution must have measure zero
support. The simplest example is the Curzon metric. This has the
potential function of a point mass. For this metric we have
\begin{gather}
  \gl = \frac{-m}{\sqrt{\gr^2+z^2}}\qquad \gm = \frac{-m\gr^2}{2(\gr^2+z^2)^2}
\end{gather}
As expected this has ADM mass of $m$. However, it is not isometric to
the Schwarzschild solution. One can calculate that the curvature
invariant $K=R_{ijkl}R^{ijkl}$ diverges as one approaches the origin
from along any straight line besides the $z$-axis. However $K$ goes to
zero along the $z$-axis.

Up to scaling there are only two value of interest for this metric,
either $m<0$ or $m>0$. The area of a level set of $\gl$
(also equal to a coordinate sphere of radius $r$) is given by
\begin{gather}
  2\pi r^2\int_0^\pi e^{\frac{2m}r-\frac{m^2\sin^2\gf}{2r^2}}\sin\gf\de{\gf}.
\end{gather}
The $m^2/r^2$ term dominates for small $r$ regardless of the sign of
$m$, so the area of these spheres goes to $0$.

The mean curvature of these spheres is given by
\begin{gather}
  -e^{\frac{2m}r-\frac{m^2\sin^2\gf}{2r^2}}\frac{m^8}{r^{18}}\left(\frac{2m^2}{r^2}
-\frac{m^3\gr}{r^5}-2(z+\gr)\right).
\end{gather}
}

The particular family of metrics we will consider are the
Zipoy--Voorhees, or $\gamma$ metrics. These are given by the potential
arising from a uniform density rod of length $2a$ and mass $m$ at the
origin. This gives a potential of
\begin{gather}
  \gl = \frac{m}{2a}\ln\frac{R_++R_--2a}{R_++R_-+2a}\qquad
\gm = -\frac{m^2}{2a^2}\ln\frac{4R_+R_-}{(R_++R_-)^2-4a^2}.
\end{gather}
Where $R_\pm = \sqrt{\gr^2+(z\pm a)^2}$. If $m=2a$, then this metric
is Schwarzschild.  It only represents the area outside the horizon.
The interval $[-a,a]$ on the $z$-axis is the event horizon.  If $m\neq
2a$ then this metric has a naked singularity at $\gr=0$, $\abs{z}\leq
a$. See \cite{ZVW-struct} for more information on the cases with
positive $m$. When $m\neq0,2a$, the resulting spacetime has ADM mass $m$
and the $[-a,a]$ interval on the $z$-axis becomes a candidate for a
negative mass singularity. To investigate the area near the
singularity, we will look at $\gr$ constant cylinders from $z=-a$ to
$a$.  The area of these cylinders is given by
\begin{gather}
2\pi\gr  \int_{-a}^{a}e^{\gm-2\gl}dz.
\end{gather}
Since both $\gl$ and $\gm$ are given by logs, the integral simplifies
to
\begin{gather}
2\pi\gr  \int_{-a}^{a}\left(\frac{4R_+R_-}{(R_++R_-)^2-4a^2} \right)^{-m^2/2a^2}
\left(\frac{R_++R_--2a}{R_++R_-+2a} \right)^{-m/a}
dz.  
\end{gather}
When $z$ is between $-a$ and $a$, and $\gr$ is small we have the
following approximations:
\begin{align}
  \frac{4R_+R_-}{(R_++R_-)^2-4a^2}&=
\frac{(a^2-z^2)^2}{a^2}\frac1{\gr^2}+\bigo(1)  
  \\
  \frac{R_++R_--2a}{R_++R_-+2a}&=
\frac{1}{4(a^2-z^2)}\gr^2+\bigo(\gr^4).  
\end{align}
Hence our integral becomes
\begin{gather}
2\pi
\gr^{m^2/a^2-2m/a+1}
4^{m/a}
\int_{-a}^a
(a^2-z^2)^{-m^2/a^2+m/a}
a^{m^2/a^2}dz.
\end{gather}
Now if $m\neq a$, the $\gr$ term has positive exponent. Thus the areas
of these surfaces go to zero. The fact that the function $e^{\gl}$ is
harmonic and goes to infinity near the $[-a,a]$ on the $z$ axis tells
us that these surfaces have zero capacity. Thus, as long as $m\neq a$,
these fulfill the definition of a negative point mass singularity. 

Continuing in this fashion we can estimate the mass of this
singularity. We will just estimate the mass of the singularities using
the level sets of the function $\gl$. Referring back to
Definition~\Ref{NPMS-mass-surf}, our function $h_i$ given by the level
set $\gl=L$ is given by
\begin{gather}
  h_i =  \frac{L}{L-1}-\frac{e^\gl}{L-1}.
\end{gather}
Here $L=\pm i$, with the sign chosen to be the opposite sign to $m$.
We can calculate $\gv(h_i)^{4/3}$ as
\begin{gather}
\gv(h_i)^{4/3}=  \frac{1}{\left(L-1\right)^{4/3}}e^{\frac43\gl
    -\frac43\gm}\left(\gl_\gr^2+\gl_z^2\right)^{2/3}.
\end{gather}
Now we will approximate the surface $\gl=L$ by a level set
$\gr=\gr_i$. Noting that $\gl_\gr$ is much larger then $\gl_z$ tells
us that this assumption is valid.
With that assumption, our mass integral becomes
\begin{gather}
  \E= \frac{2\pi\gr}{\left(e^{\gl}-1\right)^{4/3}}e^{-\frac23\gl}\int_{-a}^a
    e^{-\frac13\gm}\left(\gl_\gr^2+\gl_z^2\right)^{2/3}dz. 
\end{gather}
Now we note the first order behavior of $\gl$ and $\gm$ near $\gr=0,
\abs{z}< a$.
\begin{gather}
  \begin{aligned}
    \gl&\sim -\frac{m}{a}\ln\gr & \gm&\sim \frac{m^2}{a^2}\ln \gr\\
    \gl_\gr &\sim -\frac{m}{a}\frac1\gr & \gl_z &\sim \frac{mz}{a(z^2-a^2)}.
  \end{aligned}
\end{gather}
Using these we approximate the above integral
\begin{gather}
\E\sim \frac{2\pi\gr}{\left(\gr^{-m/a}-1\right)^{4/3}}\gr^{\frac23\frac{m}a}\int_{-a}^a
    \gr^{-\frac13\frac{m^2}{a^2}}\left(\frac{m^2}{a^2}\gr^{-2}+\frac{m^2z^2}{a^2(z^2-a^2)^2}\right)^{2/3}dz. 
\end{gather}
As $\gr$ goes to zero, this has $\gr$ dependence as
\begin{gather}
\E\sim C\cdot  \frac{\gr^{\frac23\frac{m^2}{a^2}-\frac13\frac{m}a-1}}
    {\left(\gr^{-m/a}-1\right)^{4/3}}.
\end{gather}
For some constant $C$. If $m>0$, then the bottom term contributes a
$\gr^{\frac43\frac{m}a}$ to the growth giving an overall power of
$\frac23\frac{m^2}{a^2}+\frac{m}a -1$. Otherwise it contributes
negligibly. Thus we have the following $\gr$ dependence
\begin{gather}
  \E\sim\begin{cases} C_+ \gr^{\frac23\frac{m^2}{a^2}+\frac{m}a -1} & m>0\\
  C_-\gr^{\frac23\frac{m^2}{a^2}-\frac13\frac{m}a-1} & m<0.
 \end{cases}
\end{gather}
The exponent on the $\gr$ is negative when $\frac{m}a\in
\left(-1,\frac{\sqrt{33}}4-\frac34\right)$ and negative when
$\frac{m}a$ falls outside the closure of that range. Outside that
range we have produced an example of a set of surfaces which give zero
to the mass under Definition~\Ref{NPMS-mass-surf}. Hence, since that
mass is a sup over all such surfaces, we know it must be zero.  For
$\frac{m}a$ inside that range, our surfaces give a mass of $-\infty$.
This is inconclusive, since it is entirely possible that there exists
a better behaved family of surfaces.

Checking our two known cases, $\frac{m}a=\pm 1$ we  see that the
positive Schwarzschild metric doesn't have a singularity, and the mass
of the negative Schwarzschild is nonzero and finite.  It shouldn't be
surprising that for positive $m$ outside that range we do not get a
negative mass as a Negative Point Mass Singularity since these
singularities are the only points without zero scalar curvature in a
static manifold with positive ADM mass. It seems sensible that they
shouldn't be assigned a negative mass.

\fluff{
The class of metrics we will be looking at generalizes the
Zipoy--Voorhees metrics. We will have a collection of masses
distributed along the $z$-axis, represented by a function $\gd(z)$
with compact support. The standard Zipoy--Voorhees is given by $\gd =
\frac{m}{2a}\gc_{[-a,a]}$. Given such an $\gd$ we can generate a
metric by setting
\begin{gather}
  \gl(r,z) = \int_\Re \frac{\gd(t)}{\sqrt{r^2+(z-t)^2}}\de{t}.
\end{gather}
We then generate $\gm$ by solving equations \Eqref{gm-z} and
\Eqref{gm-gr}, with the condition that $\gm\to0$ at $\infty$. 
We will call the resulting metric the Weyl--$\gd(z)$ metric. 

\section{Singularity Resolution}
For a Weyl singularity, we have a natural function to use to try and
resolve the singularity. Based on the recognition of the $\gd=-1/2$
case as negative Schwarzschild metric, we can look for a conformal
factor, $\gf$ to resolve the singularity.  By comparing the
Zipoy-Voorhees representation of Schwarzschild $m<0$ spacetime the
function $\gf=e^{-\gl}$ appears as the appropriate harmonic conformal
factor.  With this factor our model space has the spatial metric
\begin{gather}
  \tilde g=e^{2\gl}\left[\gr^2d\gth^2+
e^{2\gm}\left(d\gr^2+dz^2\right)\right].
\end{gather}

To Do:
Asymptotics of mean curvature, area and hawking mass for various
potential functions. for cylinders and \gl level sets.

}

\chapter{Open Questions}

There are a number of unanswered questions having to do with negative
point mass singularities. The most prominent is extending Theorem
\Ref{penrose-one-general} to include multiple singularities. Bray, in
\Cite{Hugh-INI}, has a solution that depends on an unproven geometric
conjecture. A further result would be what I have been calling the
``Mixed Penrose Inequality.'' This would be a result that combines the
two cases, singularities and horizons, and provides a lower bound on
the ADM mass of a manifold containing horizons and singularities. The
desired inequality is presented in Appendix~\Ref{AppB}.

It is clear that a removable singularity should have mass zero.
Precisely what conditions on an negative point mass singularity
guarantee that it is removable requires further investigation.

On a more concrete note the Zipoy--Voorhees metrics with $\frac{m}a\in
\left(-1,\frac{\sqrt{33}}4-\frac34\right)$ deserve further
investigation. In particular the metrics with $\frac{m}a$ in that
range and positive seem to offer the possibility of negative mass
singularities in a vacuum static manifold with positive ADM mass. 

On a broader scale, the current presentation of negative point mass
singularities is only in the Riemannian case. Extending the
definitions to the full Lorentzian context would require some equation
to control the behavior of the manifold in the neighborhood of the
singularity in the timelike direction. Furthermore, a new definition
of mass would have to be added, since the current definition depends
on the foliation of surfaces converging to the singularity.
Furthermore, the conformal factor definition of the regular case
doesn't translate in an obvious way to the Lorentzian case even with
the Schwarzschild metric.  Studying negative point mass singularities
in a spacetime is an important direction to pursue.

\appendix
\chapter{Miscellaneous Calculations}
\Label{append-calc}
\section{Calculation of total magnification of a NPMS lens}
\Label{append-calc-mag}
The magnification of a negative point mass singularity lens is given
by
\begin{gather}
\mu=\frac{1}{1-\frac{m^2}{\norm{x}^4}}
=\frac{\norm{x}^4}{\norm{x}^4-m^2}
=1+\frac{m^2}{\norm{x}^4-m^2}.
\end{gather}
The image locations $x_\pm$ associated to a given $y$ are
\begin{gather}
x_\pm = \frac12\left(y\pm \sqrt{y^2+4m}\right).
\end{gather}
First we calculate $\norm{x_\pm}^4$ as
\begin{align}
\norm{x_\pm}^4 &=\frac1{16}\left(y\pm \sqrt{y^2+4m}\right)^4\\
&=\frac1{16}\left(y^2\pm2y\sqrt{y^2+4m}+y^2+4m\right)^2\\
&=\frac1{16}\left(4y^4+16y^2m+16m^2+4y^2(y^2+4m)\pm(2y^2+4m)4y\sqrt{y^2+4m}\right)\\
&=\frac1{16}\left(8y^4+32y^2m+16m^2\pm\left(8y^3+16ym\right)\sqrt{y^2+4m}\right)\\
&=\frac12 y^4+2y^2m+m^2\pm\left(\frac12y^3+ym\right)\sqrt{y^2+4m}.
\end{align}
Now we define the following terms
\begin{align}
A&=\frac12y^4+2y^2m\\
B&=m^2\\
C&=y\left(\frac12y^2+m\right)\sqrt{y^2+4m}.
\end{align}
Thus $\norm{x_\pm}^4= A+B\pm C$. The negative image has negative
magnification, so we have to look at $\mu(x_+)-\mu(x_-)$. This is
given by
\begin{align}
\mu_t(x)&=\frac{\norm{x_+}^4}{\norm{x_+}^4-m^2}-\frac{\norm{x_-}^4}{\norm{x_-}^4-m^2}\\
&=\frac{A+B+C}{A+C}-\frac{A+B-C}{A-C}\\
&=1+\frac{B}{A+C}-\left(1+\frac{B}{A-C}\right)\\
&=\frac{-2BC}{A^2-C^2}.
\end{align}
We can simplify this as
\begin{align}
2BC &= ym^2\left(y^2+2m\right)\sqrt{y^2+4m}\\
A^2-C^2&=\left(\frac12y^4+2y^2m\right)^2-y^2\left(\left(\frac12y^2+m\right)\sqrt{y^2+4m}\right)^2\\
&=-y^4m^2-4y^2m^3=-y^2m^2(y^2+4m).
\end{align}
Thus we get
\begin{gather}
\mu_t= \frac{y^2+2m}{y\sqrt{y^2+4m}}.
\end{gather}

\fluff{\section{ADM mass of harmonically flat manifold}
\Label{append-calc-adm}
We want to calculate the ADM mass of a manifold that is harmonically
flat at infinity. Our manifold is given by $(\Re^3,\gf^4\gd)$. Here
$\gf$ is a harmonic function going to $1$ at infinity. Assume it has
asymptotic expansion
\begin{gather}
  \gf = 1+\frac{b}{\abs{x}}+\bigo\left(\frac{1}{\abs{x}^2}\right).
\end{gather}
We are interested in 
\begin{gather}
  m_{\ADM}=\lim_{r\to\infty}\frac1{16\pi}\int_{S_\gd^r}\left(g_{ii,j}
    -g_{ij,i}\right)n^j d\gm.
\end{gather}
Fortunately $g_{ij}=0$ if $i\neq j$ and $\gf^4$ if $i=j$. That gives
us
\begin{gather}
  m_{\ADM}=\lim_{r\to\infty}\frac{1}{2\pi}\gf^3\int_{S_\gd^r}\pp{\gf}{x^j}n^j d\gm.
\end{gather}
The gradient of $\gf$ is given by
\begin{gather}
   \grad\gf = -\frac{b}{\abs{x}^3}x+\bigo\left(\frac{1}{\abs{x}^3}\right).
\end{gather}
Looking at just the highest order terms in $\gf$ and $\grad\gf$, we get
\begin{align}
  m_{\ADM}
&=\lim_{r\to\infty}\frac{1}{2\pi}\int_{S_\gd^r}\ip{\grad\gf}{\gv}\de A\\
&=\lim_{r\to\infty}\frac{1}{2\pi}\int_{S_\gd^r}\ip{\frac{bx}{\abs{x}^3}}{\frac{x}{\abs{x}}}\de A\\
&=-\lim_{r\to\infty}\frac{1}{2\pi}\int_{S_\gd^r}
b\abs{x}^2\de A\\
&=-\lim_{r\to\infty}2b=-2b.
\end{align}}
\section{Solutions to Cusp Equation}
\Label{append-calc-cusp}
In the following we will drop the $*$ subscripts on the quantities
$Z$, $\gg$, and $m$.  We are solving the equation
\begin{gather}
   0 = Z  +\ge_\gk\left(\gg+\frac{m}{\bar z^2}\right)\bar Z,\Label{cups-eq}
\end{gather}
where 
\begin{gather}
 Z= -4i\left(\gg+\frac{m}{z^2}\right)\frac{m}{\bar z^3},
\end{gather}
and
\begin{gather}
  z=\pm\sqrt{\frac{m}{e^{-i\gf}-\gg}}.
\end{gather}
Thus we get the following
\begin{gather}
  \begin{aligned}
    \gg+\frac{m}{z^2}&=e^{-i\gf}&
    \gg+\frac{m}{\bar z^2}&=e^{i\gf}\\
    \frac{m^{3/2}}{z^3}&=\left(e^{-i\gf}-\gg\right)^{3/2}& 
\frac{\bar{m^{3/2}}}{\bar z^3}& =- \frac{ m^{3/2}}{\bar
      z^3}=-\left(e^{i\gf}-\gg\right)^{3/2}.\Label{ap-cusp-sub}
\end{aligned}
\end{gather}
Note the conjugate bar on the $m$ on the last relation. Since $m$ is
negative $m^{3/2}$ is imaginary, and so we had to multiply by $-1$
when we conjugated the previous equation. Thus our equation is
\begin{gather}
  0 =  -4i\left(\gg+\frac{m}{z^2}\right)\frac{m}{\bar z^3}
+\ge_\gk\left(\gg+\frac{m}{\bar z^2}\right)
4i\left(\gg+\frac{m}{\bar z^2}\right)\frac{m}{z^3}.
\end{gather}
Multiplying by $\frac{\sqrt m}{-4i}$ and using the substitutions in
\Eqref{ap-cusp-sub} and then multiplying by $e^{-i\gf/2}$ gives us
\begin{gather}
  0=e^{-\frac32i\gf}\left(e^{i\gf}-\gg\right)^{3/2}+\ge_\gk
  e^{\frac32i\gf}\left(e^{-i\gf}-\gg\right)^{3/2}.
\end{gather}
If we define $w= e^{-i\gf}\left(e^{i\gf}-\gg\right)$. Then we
have
\begin{gather}
  0=w^{3/2}+\ge_\gk \bar{w}^{3/2}.
\end{gather}
Which is solved when $w^{3/2}$ purely real or imaginary, when
$\ge_\gk$ is negative or positive respectively. This corresponds to
$w^3$ being purely real and having the opposite sign as $\ge_\gk$. The
imaginary part of $w^3$ is
\begin{gather}
\Iem(w^3)=\gg \sin\gf\left[4\gg^2\cos^2\gf-6\gg\cos\gf+4-\gg^2\right].
\end{gather}
Setting aside the case when $\gg=0$, we are left with the roots 
\begin{gather}
  \gf_1 = 0,\quad \gf_2=\pi, \quad \gf_{3,4} =\acos\left( \frac{3\pm
      \sqrt{4\gg_*^2-3}}{4\gg_*} \right),\quad \gf_5 =
  2\pi-\gf_3,\quad \gf_6 = 2\pi-\gf_4.
\end{gather}
If $\gg^2<3/4$ we only have at most $\gg_1$ and $\gg_2$. 
Looking at the real part of $w^3$ now
\begin{gather}
\Ree(w^3)=-4\gg^3 \cos^3\gf+6\gg^2\cos^2\gf+(3\gg^3-3\gg)\cos\gf-3\gg^2+1.
\end{gather}
Thus since there are no $\sin\gf$ terms, the real part of $w^3$ only
depends on $\cos\gf$. For each $\gf_i$, we only have a solution, and
hence a cusp, for one value of $\ge_\gk$, where $\ge_\gk$ has the same
sign as $w^3$. Thus we need to find where the sign of $w^3$ changes.
All this is in Table~\Ref{w-table}. 
\begin{table}[hptb]
\centering
\begin{tabular}{|c|c|c|c|c|}
  \hline
  $\gf_i$ & $\cos(\gf_i)$ & $\Ree(w^3)$ & $+$ & $-$\\
  \hline\hline
  $\gf_1$ & $1$ & $-\gg^3+3\gg^2-3\gg+1=(1-\gg)^3$&$\gg<1$  & $\gg>1$ \\\hline
  $\gf_2$ & $-1$ & $\gg^3+3\gg^2+3\gg+1=(1+\gg)^3$ & always & \\\hline
  $\gf_{3,5}$ & $\frac{3+\sqrt{4\gg^2-3}}{4\gg}$ &
  $1-\frac32\gg^2+\frac12\gg^2\sqrt{4\gg^2-3}$&$\gg>1$
  & $\sqrt{3/4}\leq\gg<1$  \\\hline
  $\gf_{4,6}$ & $\frac{3+\sqrt{4\gg^2+3}}{4\gg}$ &
  $1-\frac32\gg^2-\frac12\gg^2\sqrt{4\gg^2-3}$&&always \\\hline
\end{tabular}
\caption{Values of the real part of $w^3$ for various $\gf_i$.}
\Label{w-table}
\end{table}

\chapter{Work Towards the Mixed Penrose Inequality}
\Label{AppB}
\section{Purpose}
This appendix is a record of my attempts towards solving what we have
been calling the ``Mixed Penrose Inequality.''  The desired conjecture
is as follows:
\begin{conj}[Mixed Penrose Inequality]
  Let $(M, g)$ be an asymptotically flat manifold with outer
  minimizing boundary $\gS$ and negative mass singularities $\bar
  p_i$. Also assume that $(M,g)$ has nonnegative scalar
  curvature. Then
  \begin{gather}
    m_{\ADM}(G)\geq \sqrt{\frac{\abs{\gS}}{16\pi}}-\left(\sum_i
      m_{\NPMS}(p_i)^{2/3}\right)^{3/2}. 
  \end{gather}
\end{conj}
The exponent on the negative point mass singularity term comes from
consolidating the integrals for the masses as follows (assuming they
are regular):
\begin{align}
 -\left( \sum_i m_{\NPMS}[p_i]^{2/3}\right)^{3/2}
&=-\left(\sum_i\left[\frac14\left(\frac1\pi\int_{\gS_i}\gv(\gf)^{4/3}da\right)^{3/2}\right]^{2/3}\right)^{3/2}\\
&=-\left(\sum_i\frac1{4^{3/2}}\frac1\pi\int_{\gS_i}\gv(\gf)^{4/3}da\right)^{3/2}\\
&=-\frac14\left(\frac1\pi\int_{\cup_i\gS_i}\gv(\gf)^{4/3}da\right)^{3/2}.
\end{align}

This has so far been out of our reach. As a first case we have been
working on a reduced conjecture which assumes the singularity and
black hole have equal and opposite masses. Then all that is expected
is that the ADM mass is positive.
\begin{conj}
  Let $(\bar M, \bar g)$ be an asymptotically flat 3-manifold, with
  minimal boundary $\bar \gS$ and negative mass singularity $\bar p$
  of opposite mass. Also assume that $(\bar M, \bar g)$ has
  nonnegative scalar curvature. Assume that $\bar p$ can be resolved
  by a harmonic conformal factor $\gf$ so that $(M,g)$ is
  asymptotically flat, $\gS$ is still minimal, and $\gP$, the resolution
  of $p$, is minimal as well. We can assume that $\gf =1$ at
  infinity, $ \gf =0$ only on $\gP$, $\gv(\gf)=1$ on $\gP$,
  and $\gv(\gf)=0$ on $\gS$.

  Then the ADM mass of $(\bar M, \bar g)$ is non-negative. 
\end{conj}
There are three methods we have for proving Penrose style theorems.
Inverse mean curvature flow, as in \cite{IMCF}, which only takes into
account a single horizon or singularity, so probably won't be helpful
here. Bray's conformal flow of metrics as in \cite{Hugh-Penrose} is
the method we worked most with.  Lastly, Bray's Minimal Surface
Resolution technique used in \cite{Hugh-INI}, might also be useful,
however to use it, we would need to adapt it to accommodate the black
holes.  Furthermore, one of the steps in its proof is not yet
complete.

\section{Effects of Harmonic Conformal Flow}
The basic idea here is to use a harmonic conformal flow on the model
space $(M,g)$, which is equivalent to a harmonic conformal flow on the
actual space $(\bar M, \bar g)$. 

If we look at a flow on $\gf$, and set $\gps = \dd{}t \gf$, then we can
look at the changes to the various quantities. If we set $m$ to be the
ADM mass, $M$ to be the mass of the black hole, and $N$ to the mass of
the negative mass singularity.

\begin{align}
  N(\gP) &=-\frac14\left(\frac1\pi\int_\Pi
    \gv(\gf)^{4/3}dA\right)^{3/2}\\
M(\gS) & = \sqrt{\frac{\abs{\bar\gS}}{16\pi}}\left(1-\frac{1}{16\pi}\int_{\bar{
\gS}}  H^2 d\bar A\right)
\end{align}
For a first estimate we will look at the rate of change at time  zero for
initial $\gps$. First we look at $\Pi$. A direct calculation shows
\begin{align}
  \dot N(\Pi)&=-\frac 38\left(\frac1\pi\int_\Pi
    \gv(\gf)^{4/3}dA\right)^{1/2}\left(\frac1\pi \int_\Pi
  \gv(\gf)^{1/3}\left[\frac43 \gv(\gps) + H\gps\right]dA\right).
\end{align}
Under the assumption that we are working at $t=0$, we can put in the
known information about $\gf$ to get:
\begin{align}
  \dot N(\Pi)&=-\frac 32\left(\frac1\pi\int_\Pi
    \gv(\gf)^{4/3}dA\right)^{1/2}\left(\frac1\pi \int_\Pi
  \gv(\gf)^{1/3}\left[\frac43 \gv(\gps) + H\gps\right]dA\right)
\\&=N(\Pi)\frac{\int_\Pi 2 \gv(\gps) + \frac32 H\gps
    dA}{\int_\Pi dA}.
\end{align}
If $\Pi$ is minimal as well, we can drop the $H$ term. We can also
drop this term if we assume that $\gps=0$ on $\Pi$.
Likewise we can compute the change to $M$ as 
\begin{align}
  \dot M(\gS) &
  =\sqrt{\frac{1}{4\pi}}\left(\int_\gS\gf^4dA\right)^{-1/2}
  \left(\int_\gS \gps\gf^3 dA\right)\\
  & = M(\gS)\frac{\int_\gS 2\gps\gf^3 dA}{\int_\gS \gf^4 dA}.
\end{align}
The integral of $H^2$ term drops out since $H$ is zero.
We are also concerned with the change in the ADM mass of the
manifold. We can measure that with the integral of the flux across a
sphere at infinity, since this will measure the $1/r$ term of
$\gps$. (Here we assume that $\gps = 0$ at infinity.) Since the
harmonic functions have constant flux, we know that the flux across
$\gS$ and $\gP$ will also give us the $1/r$ term for $\gps$. Thus the
change in $m$ is given by
\begin{align}
  \dot m &= \frac1{4\pi}\left(\int_\gS \nabla(\gps)\cdot\gv \,dA+\int_\gP
    \nabla(\gps)\cdot\gv\, dA \right)\\
& =  \frac1{4\pi}\left(\int_\gS \gv(\gps)dA+\int_\gP
    \gv(\gps) dA \right)
\end{align}
Adding these together, the quantity we want to decrease is the
\emph{mass surplus} $X$,
\begin{align}
X &= m-M-N\\
\dot X & = \dot m -\dot M -\dot N \\
& =  \frac1{4\pi}\left(\int_\gS \gv(\gps)dA+\int_\gP
    \gv(\gps) dA \right) - \\&\qquad M\left(\frac{\int_\gS 2\gps\gf^3 dA}{\int_
\gS \gf^4 dA}\right)-N\left(\frac{\int_\Pi 2 \gv(\gps) + \frac32 H\gps
    dA}{\int_\Pi dA}\right).
\end{align}
If we can flow so that this decreases until $N$ (or $M$) is zero, then
we can apply the Penrose (NPMS) inequality to get that $X\geq 0$ at
this time, hence $X\geq0$ at $t=0$, and we are done.
At this point we tried different boundary conditions to accomplish this.

\subsection{Boundary Conditions 1}

For a first try, set $\gps$ by the boundary conditions $\gps = 0$ at
infinity, $\gps = -\frac{\ga}{2}\frac{\int\gf^4}{\int \gf^3}$ on $\gS$,
with both integrals over $\gS$, and $\gps =0$ and
$\gv(\gps)=-\frac{\gb}{2}$ on $\gP$. This gives a value for $\dot X$
of
\begin{align}
  \dot X &= \frac1{4\pi}\int_\gS \gv(\gps)dA-\frac{\gb}{8\pi}\abs{\gP}
  - \gb M-\ga N.
\end{align}

\subsection{Boundary Conditions 2}
As a second try, we can look at the boundary conditions, $\gps = 0$
at infinity, $\gps = \ga$ on $\gS$, $\gv(\gps) = \gb$ on $\gP$, where
we choose $\ga$ to be $-\frac{\int\gf^4}{2\int \gf^3}$, and pick $\gb$
so that $\dot m$ is zero. This means that the flux of $\gps$ through
the two surfaces is zero.

\subsection{Unknown Boundary Conditions}
Trying the boundary conditions $\gps = 0$ on $\gP$ and at infinity
while we leave the boundary conditions of $\gps$ on $\gS$ up in the
air for now, except that we want $\gps<0$ and $\gv(\gps)>0$.  If we
set $S$ to be a large sphere at infinity, then $\dot m$ will be given
by $\frac1{4\pi}\int_S \gv(\gps)$, with $\gv$ being the outward
pointing normal for all surfaces. Then since $\gps$ is harmonic we
know that
\begin{gather}\int_S \gv(\gps)+\int_\gS \gv(\gps)+
  \int_\gP \gv(\gps)=0\end{gather}
Now we will calculate what we get for
$\dot N$ and for $\int_\gP \gv(\gps)$ for different values of
$\gv(\gps)$. Set $\gv(\gf)=c$. Solving the formula for $N$ for
$\abs{\gP}$ we get
\begin{align}
  N &= -\frac14\left(\frac1\pi\int_\gP\gv(\gf)^{4/3}dA\right)^{3/2}\\
  & = -\frac14\frac1{\pi^{3/2}}\abs{\gP}^{3/2}\gv(\gf)^2\\
\abs{\gP} & = \left(-4\pi^{3/2}N/c^2\right)^{2/3}\\
&= 2^{4/3}\gp N^{2/3}c^{-4/3}.
\end{align}

With this in hand we can compute $\dot N$.
\begin{align}
  \dot N & =
  -\frac38\left(\frac1\pi\int_\gP\gv(\gf)^{4/3}dA\right)^{1/2}
\left(\frac1\pi\int_\gP\frac43\gv(\gf)^{1/3}\gv(\gps)dA\right)\\
  &=-\frac{\abs{\gP}^{1/2}}{2\pi^{3/2}} c
\left(\int_\gP\gv(\gps)dA\right)\\
  &=-\frac12\left(2^{2/3}\gp^{1/2} N^{1/3}c^{-2/3}\right)\frac{1}{\pi^{3/2}} c
\left(\int_\gP\gv(\gps)dA\right)\\
& = -2^{-1/3}\gp^{-1}N^{1/3}c^{1/3}\left(\int_\gP\gv(\gps)dA\right).
\end{align}
If we choose $c=-\frac{1}{32N}$. Then 
\begin{align}
  \dot N & =
  -2^{-1/3}\gp^{-1}N^{1/3}c^{1/3}\left(\int_\gP\gv(\gps)dA\right)\\
 &
 =2^{-1/3}\gp^{-1}N^{1/3}\left(\frac{1}{32N}\right)^{1/3}\left(\int_\gP\gv(\gps)dA
\right)\\
&= \frac1{4\pi} \left(\int_\gP\gv(\gps)dA\right).
\end{align}
Thus $\dot m - \dot N$ is equal to $\frac1{4\pi}\left(
  \int_\gP\gv(\gps)dA+  \int_S\gv(\gps)dA    \right)
$. Hence, in order to prove monotonicity of $X$, we only need to show
that $4\pi \dot M \geq -\int_\gS\gv(\gps)dA$.

\begin{conj}\Label{conj-2.1}
  Let $\tilde M$ be an asymptotically flat manifold with positive
  curvature. Let $\tilde M$ contain negative point mass singularities
  $p_i$ and black holes $\tilde \gS_j$. Assume that the $p_i$ can be
  resolved with a function $\gf$ on $M$ so that $\gf=1$ at infinity,
  $\gv(\gf) = 0 $ on $\gS$, $\gf=0$ on $\gP$, $\gv(\gf)= -
  \frac1{32N}$ on $\gP$, where $N$ is the mass of the $p_i$, summed
  appropriately, and $\gP$ is minimal. Then
  \begin{gather}
    m_{\ADM}(\tilde M) \geq \left(\sum_i m_{H}(\tilde
      \gS_i)^2\right)^{1/2}-\left(\sum_i
      \abs{m_{\NPMS}(p_i)}^{2/3}\right)^{3/2}.
  \end{gather}
\end{conj}
This requires only the following conjecture:
\begin{conj}
  Under the conditions of conjecture \Ref{conj-2.1}, let $\gps$ be the harmonic function
  on $M$ so that $\gps = 0$ at infinity, $\gv(\gps)=-1$ on $\gS$, and
  $\gps = 0$ on $\gP$. Then 
\begin{gather}\sqrt{4\pi}\left(\int_\gS\gf^4dA\right)^{-1/2}
  \left(\int_\gS \gps\gf^3 dA\right)\geq\int_\gS\gv(\gps)dA\end{gather}
\end{conj}

Removing this from context gives the following conjecture:
\begin{conj}
  Let $M$ be an asymptotically flat manifold with positive scalar
  curvature. Let $M$ have two sets of boundary $\gS$ and $\gP$. Let
  $\gf$ be a function so that $\gf=0$ on $\gP$, $\gv(\gf)=0$ on $\gS$,
  $\gf=1$ at $\infty$ and $\gD\gf=0$. Furthermore assume that
  $\gv(\gf)=\frac{\sqrt{\pi}}{512\abs{\gP}^{1/2}}$ on $\gP$. Then let
  $\gps$ be the harmonic function defined by $\gps=0$ on $\gP$, $\gps=0$ at
  $\infty$. We can choose our boundary conditions for $\gps$ on $\gS$ so that
  \begin{gather}0\geq\sqrt{4\pi}\left(\int_\gS\gf^4dA\right)^{-1/2}
    \left(\int_\gS \gps\gf^3 dA\right)\geq-\int_\gS\gv(\gps)dA
\end{gather}
Of course this would require that $\gv(\gps)>0$ on $\gS$ and $\gps<0$ on $\gS$.
\end{conj}

If we recast this on the original manifold $\tilde M$, we get the
following statement:

\begin{conj}
  Let $(\bar M, \bar g)$ be an asymptotically flat 3-manifold, with
  minimal boundary $\bar \gS$ and negative mass singularity $\bar p$
  of opposite mass. Also assume that $(\bar M, \bar g)$ has
  non-negative scalar curvature. Assume that $\bar p$ can be resolved
  by a harmonic conformal factor $\gf$ so that $(M,g)$ is
  asymptotically flat, $\gS$ is still minimal and $\gP$, the
  resolution of $p$, is minimal as well. We can assume that $\gf =1$
  at infinity, $\gf =0$ only on $\gP$, $\gv(\gf)=1$ on $\gP$, and
  $\gv(\gf)=0$ on $\gS$.

  Then there are boundary conditions, $X$, on $\bar \gS$ so that if 
   $\gps$ is the solution to $\bar\gD(\gps/\gf)=1$, with
  boundary values $0$ at $\infty$ and $\bar p$, and $X$ on $\bar \gS$,
  then 
  \begin{gather}0\geq\sqrt{\frac{4\pi}{\abs{\bar \gS}}}
    \left(\int_{\bar \gS} \frac{\gps}{\gf}
      dA\right)\geq-\int_{\bar\gS}\frac{\gv(\gps)}{\gf^4}dA.
\end{gather}
\end{conj}
What these $X$ boundary conditions might be we have been unable to
ascertain.

\bibliographystyle{amsplain}

\bibliography{thesis}

\providecommand{\bysame}{\leavevmode\hbox to3em{\hrulefill}\thinspace}
\providecommand{\MR}{\relax\ifhmode\unskip\space\fi MR }
\providecommand{\MRhref}[2]{%
  \href{http://www.ams.org/mathscinet-getitem?mr=#1}{#2}
}
\providecommand{\href}[2]{#2}
\begin{thebibliography}{10}

\bibitem{adm}
R.~Arnowitt, S.~Deiser, and C.~Misner, \emph{Coordinate invariance and energy
  expressions in general relativity.}, Physical Review \textbf{122} (1961),
  997--1006.

\bibitem{bartnik-mass}
Robert Bartnik, \emph{New definition of quasilocal mass}, Physical Review
  Letters \textbf{62} (1989), 2346--2349.

\bibitem{Bray-unpub}
Hubert Bray, \emph{On negative point mass singularities}, in preparation, 2007.

\bibitem{Bray-Miao}
Hubert Bray and Pengzi Miao, \emph{On the capacity of surfaces in manifolds
  with nonnegative scalar curvature}, 2007,
  \texttt{http://arxiv.org/abs/0707.3337}.

\bibitem{Hugh-Penrose}
Hubert~L. Bray, \emph{Proof of the {R}iemannian {P}enrose inequality using the
  positive mass theorem}, Journal of Differential Geometry \textbf{59} (2001),
  177--2677.

\bibitem{Hugh-INI}
\bysame, \emph{Negative point mass singularities in general relativity},
  \texttt{http://www.newton.cam.ac.uk/webseminars/pg+ws/2005/gmr/0830/bray/},
  August 2005.

\bibitem{IMCF}
G.~Huisken and T.~Ilmanen, \emph{The inverse mean curvature flow and the
  {R}iemannian {P}enrose inequality.}, Journal of Differential Geometry
  \textbf{59} (2001), 353--437.

\bibitem{ZVW-struct}
Hideo Kodama and Wataru Hikida, \emph{Global structure of the
  {Z}ippoy--{V}oorhees--{W}eyl spacetime and the $\delta=2$ {T}omimatsu-{S}ato
  spacetime}, Classical and Quantum Gravity \textbf{20} (2003), 5121--5140.

\bibitem{plw}
Arlie~O. Petters, Harold Levine, and Joachim Wambsganss, \emph{Singularity
  theory and gravitational lensing}, Birkh{\"a}user, 2001.

\bibitem{streets}
Jeffrey Streets, \emph{Quasi-local mass functionals and generalized inverse
  mean curvature flow}, Omitted from published version, 2006.

\bibitem{MTW}
Kip~S. Thorne, Charles~W. Misner, and John~Archibald Wheeler,
  \emph{Gravitation}, W. H. Freeman, 1973.

\bibitem{wald}
Robert~M. Wald, \emph{General relativity}, The University of Chicago Press,
  1984.

\end{thebibliography}

\biography
Nicholas Philip Robbins is from Merrick, New York. He received a BA
from Swarthmore College in 2001 with High Honors with the Major in
Mathematics and the Minor in Ancient History. He received a MS from
Duke University in 2003 in Mathematics. He will be a Visiting
Assistant Professor of Mathematics at St.\ Mary's College of Maryland
for 2007--2008.

\end{document}